%% file: SoergelRescalingsTakeTwo.tex
\input preamble.tex

\begin{document}

\begin{abstract} We classify gradings on the Hecke category that refine the standard $\mathbb{Z}$-grading. We also classify object-preserving autoequivalences of the Hecke category. We obtain a natural bigrading on the Hecke category which is related to the Frobenius automorphism. We also obtain an exotic grading in special characteristic that can be used to categorify many Hecke algebras with unequal parameters, including all Hecke algebras with unequal parameters for all finite and affine Weyl groups. This paper is a replacement for \cite{EBigraded}, which is now obsolete.
 \end{abstract}

\title{Gradings on the Hecke category, and categorification with unequal parameters}

\author{Jonas Antor} \address{University of Bonn, Germany}

\author{Ben Elias} \address{University of Oregon, Eugene}

\maketitle

\setcounter{tocdepth}{1}
\tableofcontents

\section{Introduction}

\subsection{Overview}

The goal of this paper is to explore (multi-)gradings on the Hecke category, which is a $\Z$-graded monoidal categorification of the Iwahori-Hecke algebra. The Hecke category depends on a choice of realization $V$ (something like a reflection representation with a choice of simple roots) of a Coxeter system $(W,S)$. Each $s \in S$ is assigned a color in the diagrammatic Hecke category, so we refer to colors and choices of $s \in S$ interchangeably.

We classify all possible gradings on the Hecke category which refine the original $\Z$-grading. We prove that all such gradings are determined by the grading on the
univalent vertices (certain generating morphisms), and the grading on $V$, which must be $W$-invariant. Consequently there are three routes by which one could generalize the original
$\Z$-grading. The first route is to break color symmetry, assigning different degrees to univalent vertices based on their color. For example, let $\CB$ denote the set of connected components of the Coxeter graph. It is fairly obvious
that one can assign different degrees to different connected components in the Coxeter graph\footnote{There are exceptions even to this, see Example \ref{ex:yuck}.}, yielding a grading by $\Z^{\CB}$. Beyond that, it is
subtle to determine which colors can have independent degrees, and which colors must be linked. The second route is to break the vertical (and rotational) symmetry of the grading,
assigning different degrees to the ``startdot'' and ``enddot'', the two different versions of a univalent vertex. The third route is to choose a refined $W$-invariant grading on $V$.
As the grading on simple roots is determined by the grading on univalent vertices, and as $V$ is typically generated by simple roots and $W$-invariant vectors, this generalization is
rather dull, merely allowing one to assign new degrees to $W$-invariant vectors. Of course, all three generalizations can be combined.

Beyond this classification, here are the results of the paper in summary, where for simplicity we ignore the grading on $V$. A grading is \emph{inner} over another grading if it can be obtained by merely renaming (i.e. shifting) the generating objects. \begin{itemize}
\item We compute the Grothendieck group when equipped with a more general grading.
\item For most realizations, every grading is inner over the grading by $\Z^{\CB}$.
\item For heretofore unexplored realizations in finite characteristic, there are additional gradings. We use these to categorify Hecke algebras with unequal parameters.
\item Breaking vertical symmetry gives rise to a bigrading which is related to the action of Frobenius on morphisms between perverse sheaves.
\item Every object-fixing autoequivalence of the Hecke category comes from a grading.
\end{itemize}
This last point gives an easy criterion for when two such equivalences are equal. In \cite{EHDrinfeld} this is used to simplify the computation that conjugation with the half twist comes from the bigrading just mentioned.

\begin{remark} The bigrading and observations about object-fixing autoequivalences were originally found in \cite{EBigraded}. This paper replaces \cite{EBigraded}, which did not properly observe the potential for breaking color symmetry and categorifying unequal parameter Hecke algebras. \end{remark}

%

Henceforth we assume the reader is familiar with the diagrammatic Hecke category $\HC$ associated to a realization $V$ of a Coxeter system $(W,S)$, over a commutative base ring
$\Bbbk$. An explicit presentation can be found in \cite[\S 10.2.4]{EMTW} or \cite[\S 5.1]{EWGr4sb}. Typically, this is defined as a $\Z$-graded category, where e.g. univalent vertices
have degree $1$, and elements of $V$ have degree $2$ when viewed within the polynomial ring $R = \Sym(V)$. Note that $R \cong \End(\1)$, where $\1$ is the monoidal identity. We make
use of the Demazure operators $\pa_s$ for each $s \in S$, see \cite[\S 4.3]{EMTW} or \cite[\S 3.3]{EWGr4sb}.

In the introduction we assume that $(W,S)$ has no parabolic subgroup of type $H_3$. The diagrammatic Hecke category $\HC$ is not currently defined in types containing $H_3$, because
the appropriate Zamolodchikov relation is unknown. We do prove that, when this relation is found, it will be compatible with the new gradings in this paper.

{\bf Acknowledgments} We wish to thank Daniel Tubbenhauer, Matt Hogancamp, Olivier Dudas and Paul Wedrich. The first author was supported by NSF grants DMS-2201387 and DMS-2039316.

\subsection{Breaking vertical symmetry}

The Hecke category can be bigraded, by breaking the vertical symmetry. The second grading is in some sense the Soergel-ification of the Forbenius action on morphisms of constructible sheaves, as we explain in \S\ref{sec:frobenius}.

\begin{thm} \label{thm:bigraded} The diagrammatic Hecke category $\HC$ (if it is well-defined, see below) can be defined as a $\Z^2$-graded monoidal category. The bidegrees of the generating morphisms are as follows.
\begin{subequations} \label{eq:bigrading}
\begin{equation} \label{eq:onecolorgens} \deg \left( \startdotred \right) = (1,0), \qquad \deg \left( \finaldotred \right) = (0,1), \qquad \deg \left( \splitred \right) = (0,-1), \qquad 
	\deg \left( \mergered \right) = (-1,0), \end{equation}
\begin{equation} \label{eq:degofmst} \deg \left( \sbotttop \right) = (0,0). \end{equation}
\end{subequations}
The picture in \eqref{eq:degofmst} represents an arbitrary $2m_{st}$-valent vertex. The degree of any element of $V$, viewed as a linear polynomial in $R$, is $(1,1)$. The original degree function is obtained from the bidegree via the map
\begin{equation} \Z^2 \to \Z, \qquad (m,n) \mapsto m+n. \end{equation}
\end{thm}

The startdot is the first generator depicted in \eqref{eq:onecolorgens}, and the enddot is the second. There is one such generator for each color $s \in S$.

\begin{rem} Note that cups and caps have bidegree $(1,-1)$ and $(-1,1)$ respectively. The bidegree of a morphism is an isotopy invariant relative to a fixed boundary, but the bidegree is not preserved under adjunction isomorphisms.  \end{rem}


\begin{rem} This bigrading extends to the 2-category known as the singular (diagrammatic) Hecke category, which describes singular Soergel bimodules, see \cite[\S 24]{EMTW}. The
idea that the singular Hecke category could be bigraded originated out of \cite[\S 3.4]{EKo}, where two separate length functions $\ell^-$ and $\ell^+$ are defined on double
cosets. \end{rem}

\subsection{Breaking color symmetry} \label{subsec:colorsym}

\begin{notation} Let $\Z^S$ denote the free abelian group with a basis $\{1_s\}_{s \in S}$. If $\sim$ is an equivalence relation on $S$, let $\Z^{S/\sim}$ denote the free abelian group
on equivalence classes, which is also the quotient of $\Z^S$ by the relations $1_s = 1_t$ for $s \sim t$. By abuse of notation, we still let $1_s$ denote the image of $1_s$ in
$\Z^{S/\sim}$. \end{notation}

\begin{question} \label{question:color} For which equivalence classes $\sim$ is there a grading on $\HC$ valued in $\Z^{S/\sim}$, where an $s$-colored dot has degree $1_s$, an
$s$-colored trivalent vertex has degree $-1_s$, and $2m$-valent vertices have degree zero? (We also assume the grading on $\HC$ restricts to a $W$-invariant grading on $V$.) \end{question}

The original grading on $\HC$ is an example of such a grading, for the trivial equivalence relation where everything is equivalent. To break color symmetry would
be to find such a grading for a more refined equivalence relation. The following simple lemma governs this question, though the implications are subtle.

\begin{lemma} \label{lem:fdegree} In the setting of Question \ref{question:color}, let $f \in V$ be homogeneous such that $s(f) \ne f$. Then $\deg(f) = \deg(\alpha_s) = 2_s := 1_s + 1_s$. \end{lemma}

\begin{proof} The $s$-colored barbell (a composition of dots) is equal to $\alpha_s$ and has degree $2_s$. Recall that $f-s(f) = \pa_s(f) \alpha_s$. Since both sides are homogeneous,
either they are both zero or they have the same degree. \end{proof}
	
Recall that $s(\alpha_t) = \alpha_t + [2]_{s,t} \alpha_s$ for a scalar $[2]_{s,t} := -\pa_s(\alpha_t) = -\langle \alpha_s^\vee, \alpha_t \rangle \in \Bbbk$ which is found (up to sign) in the Cartan
matrix of the realization. If $[2]_{s,t} \ne 0$ or $[2]_{t,s} \ne 0$, then $\alpha_s$ and $\alpha_t$ must have the same degree.

\begin{defn} Fix a realization $V$ of a Coxeter system $(W,S)$. Let $\sim_V$ be the equivalence relation on $S$ generated by
\begin{equation} \label{eq:simdef1}
    s \sim_V t \text{ if } [2]_{s,t} \ne 0 \text{ or } [2]_{t,s} \ne 0.
\end{equation}
Note that $\sim_V$ depends only on the Cartan matrix of the realization $V$, not on $V$ itself. We write $\VB$ rather than $S/\sim_V$ for the equivalence classes of $\sim_V$.
\end{defn}

To summarize, $\sim_V$ is the coarsest equivalence relation such that $\HC$ could conceivably be graded by $\Z^{S/\sim}$, since $s \sim_V t$ implies $\deg(\alpha_s) =
\deg(\alpha_t)$ implies $2_s = 2_t$. The question of whether a grading by $\Z^\VB$ is attainable is more subtle and depends on $V$ not just its Cartan matrix; this topic discussed in \S\ref{subsec:polyssuck}. To avoid this technical discussion, in the introduction we assume that $V$ is spanned by the simple roots. In Theorem \ref{thm:mostgeneral} we prove that, whenever it is well-defined, $\HC$ can indeed be graded by $\Z^\VB$.

\begin{notation} Let $\sim_C$ denote the equivalence relation on $S$ generated by $s \sim_C t$ if $m_{st} \ne 2$. Equivalence classes for $\sim_C$ are called \emph{connected
components}. We write $\CB$ rather than $S/\sim_C$ for the set of connected components. \end{notation}

By definition of a realization, if $m_{st} = 2$ then $[2]_{s,t} = [2]_{t,s} = 0$. Thus $s \not\sim_C t$ implies $s \not\sim_V t$, so each equivalence class for $\sim_V$ is contained in
a single connected component. Note also that $[2]_{s,t} = [2]_{t,s} = 0$ implies that $s$ and $t$ commute, so $m_{st}$ cannot be odd, and the realization is not faithful if $m_{st} >
2$. For ``typical'' realizations, $\sim_V$ and $\sim_C$ are the same relation! Typical here is used in the colloquial sense, but includes all faithful
realizations, and even non-faithful, characteristic $p$ specializations of the usual integral realizations of crystallographic groups. After all, in a standard Cartan matrix, either
$[2]_{s,t} = 1$ or $[2]_{t,s} = 1$ whenever $m_{st} \in \{3,4,6\}$.

It was already relatively straightforward to generalize the original grading to one valued in $\Z^{\CB}$, called the \emph{component grading}. The Hecke category is a tensor product
(over the polynomial ring) of the Hecke categories associated to each connected component, and the component grading is merely the product of the ordinary grading for each component.
For typical realizations, this obvious generalization is the best one can do to break color symmetry.

In Theorem \ref{thm:mostgeneral} we construct atypical realizations where $\sim_V$ strictly refines $\sim_C$. We introduce this in \S\ref{subsec:atypical}.

\subsection{Classification} \label{subsec:combine}

One might attempt to combine the vertical and color symmetry breaking with a multigrading valued in \[ \prod_{s \in S} \Z^2 = (\Z^2)^S,\] where the $s$-colored startdot is
assigned grading $(1,0)$ in the $s$-th factor, denoted as $(1,0)_s$. Clearly this will not work as stated: as in the previous section, we require $\deg(\alpha_s) = \deg(\alpha_t)$ when $s \sim_V t$, meaning that we need $(1,1)_s =
(1,1)_t$.

\begin{defn} \label{lambdadef} Let $\sim$ be an equivalence relation on $S$. Let $\Lambda_{\sim}$ be the quotient of $(\Z^2)^S$ by the relation that $(1,1)_s = (1,1)_t$ whenever $s \sim t$. Let
$(\Z^2)^{S/\sim}$ be the further quotient of $(\Z^2)^S$ by the relations that $(1,0)_s = (1,0)_t$ and $(0,1)_s = (0,1)_t$ whenever $s \sim t$. We write $\Lambda_{\CB}$ instead of $\Lambda_{\sim_C}$, and $\Lambda_{\VB}$ instead of $\Lambda_{\sim_V}$.
\end{defn}

When $V$ is spanned by the simple roots, then $\HC$ can be graded by $\Lambda_{\VB}$. For a further generalization when $V$ is not spanned by simple roots, see Theorem \ref{thm:
general grading by components}. We prove that this is the most general possible grading refining the original $\Z$-grading. Again, recall that typically $\VB = \CB$.

While the grading by $\Lambda_{\CB}$ is new, it is also, in some sense, boring. Suppose one had an ungraded category, and then equipped it with a grading by an abelian group $A$ so that each
pair of objects $X$ and $Y$, the morphism space $\Hom(X,Y)$ is homogeneous in a single degree $d(X,Y) \in A$. We call this an \emph{inner grading}. Such a grading contains essentially no information about morphisms and does not
produce any nontrivial decompositions of morphism spaces. In fact, the grading contains only information about objects: one can associate\footnote{This is not quite true, see
Remark \ref{rmk:maybeneedtoextendgroup}.} a degree in $A$ to each object, such that $d(X,Y) = \deg(Y) - \deg(X)$.  Inner gradings do not have a significant effect on the
Grothendieck group: they merely add a new formal grading parameters, and rescale the meaning of each object accordingly. The concept of inner gradings extends to monoidal categories
under the requirement that $\deg(X \ot Y) = \deg(X) + \deg(Y)$.

One can also discuss whether a grading is inner relative to an original grading that it refines. In our context, an $A$-grading on a graded category would be relatively inner if each
original homogeneous piece $\Hom^k(X,Y)$ is now placed in a single degree $d(k,X,Y) = \deg(Y) - \deg(X) + \deg(k) \in A$. The bigrading on the Hecke category is indeed inner, relative to
the original grading. The grading by $\Lambda_{\CB}$ is inner relative to the component grading by $\Z^{\CB}$. Our classification of gradings effectively concludes that breaking
vertical symmetry will only produce inner gradings, a ``no surprises'' result.

\begin{remark} For the Hecke category it is desirable that the object $B_s$ remain self-dual. Shifting the grading on $B_s$ might violate self-duality. However, it is possible that the ``bar involution'' extends to the new grading group in such a way that $\deg(B_s)$ is fixed, in which case $B_s$ will still be self-dual with the new inner grading. This is the case for the bigrading in Theorem \ref{thm:bigraded}. We discuss this further in \S\ref{subsec:rescaled}. \end{remark}

We pair these grading classification results with a classification of object-preserving autoequivalences of $\HC$ in \S\ref{sec-autoequiv}. We prove that every such autoequivalence
will rescale homogeneous morphisms based on their degree in some grading. To verify that two object-preserving autoequivalences are the same, one need only check that they agree on each startdot and on $V
\subset \End(\1)$. For an example application, conjugation by the full-twist Rouquier complex is equivalent to some object-preserving autoequivalence
of $\HC$, and this result dramatically streamlines the work needed to prove that it is the identity functor, see \cite{EHDrinfeld}.

To summarize, for typical realizations the most general grading on $\HC$ will be an inner extension of the component grading, just as the bigrading of Theorem \ref{thm:bigraded} is an
inner extension of the original grading. However, just because the bigrading is inner does not mean it is not interesting or useful, as there are connections between the bigrading and
the Frobenius operator and with conjugation by the half twist. Meanwhile, the classification result itself, that every grading is boring, can help one quickly pin down gradings or
autoequivalences, which has useful applications.

For most users of the Hecke category, the results discussed above and their applications will be the main utility of this paper. The proofs are relatively straightforward as well.
However, for those readers willing to entertain unusual realizations, there are surprises! There are yet more possible gradings, with significant consequences. The proof of the main
result below requires new techniques as well.

\subsection{Jones Wenzl projectors at zero, and unlinking colors} \label{subsec:atypical}

Let us return to the setting of Question \ref{question:color}. We have already noted that $\sim$ must be coarser than $\sim_V$, i.e. that if $[2]_{s,t} \ne 0$ then $s \sim t$. But is
this the only condition we need? There are subtleties with polynomials, see \S\ref{subsec:polyssuck}, so let us ignore them by assuming $V$ is spanned by simple roots. In this case can
we grade $\HC$ by $\Z^{\VB}$? The answer is yes, and it is not obvious. The same discussion applies to a more general grading by $\Lambda_\VB$, where the result is stated as Theorem
\ref{thm:mostgeneral}.

To explain this theorem, we first recall the conditions under which $\HC$ is well-defined, as outlined in \cite[\S 5.1]{EWLoc}. First, the space $V$ must be a representation of $W$,
which when $[2]_{s,t} = [2]_{t,s} = 0$ is equivalent to $m_{st}$ being either even or infinite. Second, for each pair $s, t \in S$ with $m_{st} < \infty$ we require the existence and
rotatability of the two-colored Jones-Wenzl projector $\JW_{m_{st} - 1}$. A precise condition for existence and rotatability was given by Hazi in \cite{HaziRotatable}, in terms of properties of
two-colored quantum numbers like $[2]_{s,t}$, and their associated two-colored quantum binomial coefficients. When $[2]_{s,t} = [2]_{t,s} = 0$, Hazi's condition simplifies. The
following lemma is proven as Lemma \ref{lemma: m is 2 p^k}.

\begin{lemma} Suppose that $\HC$ is well-defined, and $s, t \in S$ are such that $[2]_{s,t} = [2]_{t,s} = 0$. Then either $m_{st} \in \{2,\infty\}$, or $\Bbbk$ has finite characteristic $p$ and $m_{st} = 2p^k$ for some $k \ge 1$. \end{lemma}

%

When $m_{st} = 2p^k$ and $[2]_{s,t} = [2]_{t,s} = 0$ we \emph{can} break color symmetry! The hardest part of the proof is showing the homogeneity of the relation\footnote{In the
body of the paper we instead study the similar relation which resolves a dotted $2m_{st}$-valent vertex.} which rewrites a composition of two $2m_{st}$-valent vertices as a linear combination of other
diagrams. These other diagrams are built from dots and trivalent vertices, and arise as deformation retracts of two-colored crossingless matchings. In the usual $\Z$-grading, all such deformation
retracts have degree zero. After breaking color symmetry, the retracts of different crossingless matchings have different gradings! The relation states that a doubled $2m_{st}$-valent vertex is a
retract of $\JW_{m_{st}-1}$. The crucial new observation is that each relevant crossingless matching either has degree zero, or has coefficient zero in $\JW_{m_{st}-1}$.

\begin{rem} When studying the Temperley-Lieb algebra with parameter $[2] = q + q^{-1}$ and its Jones-Wenzl projectors, there are two non-semisimple specializations which appear most often in the literature: the one where $q=\pm 1$ and $\Bbbk$ has characteristic $p$, and the one where $\Bbbk$ has characteristic zero but $q$ is a primitive $2k$-th root of unity. In the former $[p]=0$ and in the latter $[k]=0$, which is their first similarity. When $\Bbbk$ has finite characteristic and $q$ is a nontrivial root of
unity one is in so-called \emph{mixed characteristic}; see e.g. \cite{SuttonTiltMod} for a study of Jones-Wenzl projectors in this setting. This is the relevant
setting for us, as $\Bbbk$ has characteristic $p$, and $[2]=0$ correponds to $q$ being a fourth root of unity. \end{rem}

Now fix a Coxeter system $(W,S)$ and a base field $\Bbbk$ with characteristic $p$, possibly $p=0$. There is a realization $V$ of $(W,S)$ over $\Bbbk$ where $\HC$ is well-defined such
that $[2]_{s,t} = [2]_{t,s} = 0$ whenever this is possible (see Lemma \ref{lem:coloradapted}), so that the grading can be generalized as much as possible under the constraints
discussed above. In other words, there is a realization where the equivalence relation $\sim_V$ (which depends on the realization) agrees with the equivalence relation $\sim_{\Bbbk}$
below. Note that this will \emph{not} be the finite characteristic specialization of the usual Cartan matrix for crystallographic groups.

\begin{defn} \label{defn:simp} Fix a Coxeter system $(W,S)$ and a field $\Bbbk$. Suppose that $\Bbbk$ has characteristic $p$, possibly zero. If $p > 0$ let $\sim_{\Bbbk} = \sim_p$ be the equivalence relation on $S$ generated by $s \sim_p t$ if $m_{st} < \infty$ and $m_{st} \neq 2p^k$ for all $k$. If $p=0$ let $\sim_{\Bbbk} = \sim_0$ be the equivalence relation on $S$ generated by $s \sim_0 t$ if $m_{st} \notin \{2,\infty\}$.
\end{defn}

\begin{rem} If desired one can modify the definition of $\sim_{\Bbbk}$ so that $m_{st} = \infty$ is treated separately (e.g. $m_{st} = \infty$ implies $s \sim_{\Bbbk} t$), and one can
still construct a realization $V$ over $\Bbbk$ such that $\sim_V = \sim_{\Bbbk}$. \end{rem}

\subsection{Computing the Grothendieck group} \label{intro-grothgroup}

Let $A$ be an abelian group and let $\Z[A]$ denote the group algebra of $A$. The (split) Grothendieck group of an $A$-graded additive Karoubian category is naturally a
$\Z[A]$-algebra, with the action of $A$ coming from grading shifts. For this reason the ordinary grading on the Hecke category equips its Grothendieck group (the Hecke algebra)
with an action of $\Z[v,v^{-1}] = \Z[\Z]$. Regardless of the realization, this Grothendieck group is free of rank $\# W$ over $\Z[v,v^{-1}]$, which is a consequence of the
classification of indecomposable objects in $\HC$. More generally for an $A$-grading we can ask: is the Grothendieck group necessarily a free $\Z[A]$-algebra of rank $\# W$? How
are direct sum decompositions in $\HC$ affected by the choice of grading?

We are \emph{not} only interested in gradings which refine the original grading! For example, consider the specialization of $Z^{S/\sim}$ to $\Z$ which sends $1_s \mapsto 1$ and $1_t \mapsto 2$ for some $s \ne t$.

When a grading is inner relative to the original grading, one can compute its Grothendieck group easily: it is isomorphic to the original after base change, but the meaning of each generator is rescaled. Our computation of Grothendieck groups is general enough that it need not rely on this trick, but we use it for pedagogical reasons.

While any homogeneous idempotent in an $A$-graded category must have degree zero, it is possible for an idempotent to be non-homogeneous. A simple example of this phenomenon is a
non-diagonal idempotent $2 \times 2$ matrix (using the usual grading by diagonal on the space of matrices). Even assuming one has enough decompositions by homogeneous idempotents,
when computing Grothendieck groups it is important to consider \emph{factored} idempotent decompositions, where idempotents are expressed as a composition of projection and
inclusion maps. Factorization can be used to argue that the image of an idempotent is actually isomorphic to another, existing object. In theory, factored decompositions can be
affected by the grading, as projection and inclusion maps need not be homogeneous.

The original categorification results (see \cite[Theorem 6.25]{EWGr4sb} or \cite[Theorem 11.39]{EMTW}) for the diagrammatic Hecke category were proven using the double leaves basis and
the technology of cellular categories. Rather than try to bootstrap these results to deduce something about the new gradings, we instead prove the results for the new gradings ab ovo
via the same technology. Let us assume that there is a group homomorphism $A \to \Z$ for which the degree of all univalent vertices is sent to a positive number. Under this hypothesis,
in Theorem \ref{thm:categorification theorem} we compute the Grothendieck group of the $A$-graded Hecke category, and prove that it is a free $\Z[A]$-algebra of the expected rank. From this, we deduce that the Grothendieck group is isomorphic to a Hecke algebra with unequal parameters.

\subsection{The Hecke algebra of (affine) Weyl groups}

A case of particular interest is when $W$ is a finite or affine Weyl group. The following notion of 'special prime' is borrowed from terminology for the corresponding reductive group (c.f. \cite{SteinbergRep}).
\begin{defn} The prime $p=2$ is \emph{special} when $W$ has finite or affine type $B_n, C_n, F_4$. The prime $p=3$ is \emph{special} for finite or affine type $G_2$.
\end{defn}
When $p$ is special for $W$, we can find a realization for which $[2]_{s,t} = [2]_{t,s} = 0$ whenever $m_{st} > 2$. The following theorem is proved as Corollary \ref{corollary: categorification result for p grading}.

\begin{theorem}
Suppose that $\Bbbk$ has characteristic $p$ (possibly $0$). Then there is a realization $V$ of $W$ and a $\mathbb{Z}^{S/\sim_p}$-grading on $V$ and $\HC$ such that $\HC$ categorifies the Hecke algebra of $W$ with generic unequal parameters $\{ v_s \mid s \in S \} $ subject to the relation $v_s = v_t$ whenever $s \sim_p t$. In particular, if $W$ is an (affine) Weyl group and $p$ is special for $W$, this yields a categorification of the Hecke algebra of $W$ with arbitrary generic unequal parameters.
\end{theorem}

By specializing the grading above one can get non-generic unequal parameters.

It would be interesting to compute the basis of the Grothendieck group corresponding to the indecomposble objects in $\HC$. We call this the \emph{double-0-$p$ canonical basis}\footnote{The name is inspired by the special case $p=7$.}, because it is associated to a Cartan matrix where both off-diagonal entries are zero\footnote{In contrast, the $p$-canonical basis is associated to the usual Cartan matrix, where for special primes only one off-diagonal entry is zero.}. For the infinite dihedral group one can consider the same question in characteristic zero, and we compute its \emph{double-0 canonical basis} in \S\ref{section:dihedral example}. In particular, this example shows that the new basis we obtain need not agree with the Kazhdan-Lusztig basis for unequal parameters. This is to be expected since the Kazhdan-Lusztig basis with unequal parameters does not satisfy positivity.

The categorification result above can also be thought of as a Soergel bimodule analogue of the geometric realization of unequal parameter affine Hecke algebras from \cite{AntorGeoRel}. Interestingly, the geometric set-up from \cite{AntorGeoRel} also requires special characteristic $p$, though for a different reason: it exploits that the adjoint representation of a reductive group becomes reducible in special characteristic. It would be interesting to see if there is a categorical equivalence that relates the two constructions, similar to the Bezrukavnikov equivalence for the equal parameter setting \cite{BezruTwoGeoRel}.

\section{Gradings on the Hecke category}

\subsection{Reduction to the Jones-Wenzl relations}

\begin{defn} We say that a grading\footnote{What is an $A$-graded category? There are several standard conventions which are easy to pass between, see e.g. \cite[\S 11.2.1]{EMTW}. For simplicity we use the convention that morphism spaces are $A$-graded, and one can shift objects by $A$, with the obvious compatibility between these structures.} of $\HC$ by an abelian group $A$ \emph{refines the original grading} if there is a group homomorphism $\For \co A \to \Z$, and for any morphism $f$ which is homogeneous of degree $a$ in the $A$-grading, $f$ is also homogeneous of degree $\For(a)$ in the original $\Z$-grading. \end{defn}

Our goal is to find all possible multigradings on $\HC$ which refine the original grading (and where the base ring $\Bbbk$ has degree zero). The following proposition essentially reduces this question to checking whether the Jones-Wenzl relation is homogeneous.

\begin{proposition}\label{prop: gradings extend if and only JW is homogeneous}
    \begin{enumerate}
        \item Assume that there is a grading on $\HC$ valued in an abelian group $A$, which refines the original $\mathbb{Z}$-grading. Then this restricts to a $W$-invariant $A$-grading on $V \subset R = \End(\1)$, and we can find $f_s, g_s \in A$ for each $s \in S$ such that
        \begin{subequations}  \label{eq:multigrading}
        \begin{equation} \label{eq:onecolor} \deg \left( \startdotred \right) = f_s, \qquad \deg \left( \finaldotred \right) = g_s, \qquad \deg \left( \splitred \right) = -g_s, \qquad 
	   \deg \left( \mergered \right) = -f_s, \end{equation}
        \begin{equation} \label{eq:2mvalent} \deg \left( \sbotttop \right) = \begin{cases} 0 & \text{ if } m_{st} \text{ is even}, \\ g_s - g_t & \text{ if } m_{st} \text{ is odd}. \end{cases} \end{equation}
        \end{subequations}
        Moreover, for $s \ne t \in S$ we have
\begin{subequations} \label{eq:polygradingcriterion}
	\begin{equation} \alpha_s \text{ is homogeneous, } \qquad \deg(\alpha_s) =f_s + g_s, \end{equation}
	\begin{equation} f_s +g_s = f_t + g_t \text{ whenever } \langle \alpha_s, \alpha_t^{\vee} \rangle\neq 0. \end{equation}
\end{subequations}
    \item Conversely, assume we are given an abelian group $A$ and a $W$-invariant $A$-grading on $V$, together with elements $g_s, f_s \in A$ such that \eqref{eq:polygradingcriterion} holds. Then this extends to an $A$-grading on $\HC$ as in \eqref{eq:multigrading} if and only if all Jones-Wenzl relations are homogeneous.
    \end{enumerate}
\end{proposition}
\begin{proof}
We first prove (1). Suppose we have an $A$-grading on $\HC$ which refines the original grading. Since $\startdotred$ lives in a one-dimensional graded morphism space for the original grading, it must be homogeneous for the $A$-grading. Write $f_s$ for its new degree. By the same argument, each non-polynomial generator of $\HC$ is homogeneous for the $A$-grading. Write $g_s$ for the new degree of $\finaldotred$. Our arguments below do not use the assumption that the $A$-grading refines the original grading, only the consequence that the non-polynomial generators are homogeneous, and that $V$ is a sum of its homogeneous components.

Consider the unit and counit relations \cite[(8.5ab)]{EMTW}. Each side of the relation is a diagram which is now known to be homogeneous for the new grading, but the degrees of the two sides might be unequal. Note that it is a priori possible for $\HC$ to have a grading where the relations are not homogeneous! However, each homogeneous component of the relation must also be a valid equality in the category. If the gradings on either side of \cite[(8.5ab)]{EMTW} disagree, then both sides are individually zero, and we must have that $s$-colored strand (an identity map) is zero. This is not true in $\HC$ (we are not permitted to impose additional relations, which would be studying gradings on a quotient of $\HC$ and not on $\HC$ itself). Thus the unit and counit relations must be homogeneous, from which we compute that the degree of the split must be $-g_s$, and the degree of the merge must be $-f_s$, as in \eqref{eq:onecolor}. In similar fashion, consider the barbell relation \cite[(10.7c)]{EMTW}. If $\alpha_s$ were not homogeneous of degree $f_s + g_s$, then by projecting this relation to some other homogeneous component, some nonzero polynomial would be zero, a contradiction. Hence $\alpha_s$ is homogeneous of degree $f_s + g_s$.

Now consider each relation in the presentation of $\HC$, listed in \cite[(10.7)]{EMTW} and \cite[(10.14)]{EMTW}. With the exception of \cite[(10.7d)]{EMTW}, each relation is a
linear combination of diagrams which are now known to be homogeneous, but not necessarily known to have the same grading. In each case except one (by examination), if the gradings
were distinct then a homogeneous component of the relation would be an equality which does not hold in $\HC$, giving a contradiction. This is even true of the Jones-Wenzl relation
(it is a linear combination of diagrams, and spans the one-dimensional space of relations between those diagrams). The exception is the Zamolodchikov relation of type $H_3$, which
is not known explicitly.

Suppose that $f \in R$ is homogeneous. The argument that \cite[(10.7d)]{EMTW} is homogeneous for the new grading is similar but requires slightly more thought. Rewrite $sf$ as $f - \pa_s(f) \alpha_s$. Rewrite $\pa_s(f) = g + h$, where $h$ is its homogeneous component in degree $\deg(f) - \deg(\alpha_s)$ (possibly $h =0$), and $g$ is the sum of the remaining homogeneous components. By subtracting the homogeneous component of \cite[(10.7d)]{EMTW} in degree $\deg(f)$ from the original relation, we get a two-term relation: $g \alpha_s$ to the right of the identity of $B_s$, plus $g$ next to a broken line (two dots), is zero. However, the identity and the broken line are linearly independent for the right $R$-action in $\HC$, whence we deduce that $g=0$. Now it is simple to conclude that \cite[(10.7d)]{EMTW} is homogeneous. This implies that $s(f)$ is also homogeneous, and $\deg(s(f)) = \deg(f)$.

Since $V$ is homogeneous with respect to the original grading, it must split into homogeneous components with respect to the new grading. The previous paragraph proves that
$\deg(s(v)) = \deg(v)$ for any $v \in V$ homogeneous. Consequently, the new grading on $V$ is $W$-invariant. If $v \in V$ is homogeneous and $0 \ne \pa_s(v) \in \Bbbk$, then $v - s(v)$ is homogeneous and $v - s(v) = \pa_s(v) \alpha_s$ is nonzero, so $\deg(v) = \deg(\alpha_s)$. Pick $s \ne t$ and recall that $-[2]_{s,t} = \langle \alpha_t, \alpha_s^{\vee} \rangle = \pa_s(\alpha_t)$. Thus if $[2]_{s,t} \ne 0$ then $\deg(\alpha_s) = \deg(\alpha_t)$. So \eqref{eq:polygradingcriterion} holds.


Pick $s \ne t$ with $m_{st} < \infty$. The $2m_{st}$-valent vertex also lives in a one-dimensional graded morphism space, so it is homogeneous for the new grading. We set
\begin{equation} \deg \left( \sbotttop \right) = h_{s,t}. \end{equation}
A priori, it is possible that $h_{s,t} \ne h_{t,s}$.

When $m_{st}$ is even, the two-color associativity relation \cite[(10.7fh)]{EMTW} is homogeneous if and only if $h_{s,t} - g_s = 2 h_{s,t} - g_s$, or equivalently, $h_{s,t} = 0$. When $m_{st}$ is odd, the two-color associativity relation \cite[(10.7g)]{EMTW} is homogeneous if and only if $h_{t,s} - g_s = 2h_{t,s} - g_t$, or equivalently, $h_{t,s} = g_t - g_s$. Thus the degree $h_{s,t}$ agrees with \eqref{eq:2mvalent}. We have now proven (1).

Now we prove (2). We have already argued above that all relations (except the Zamolodchikov relation in type $H_3$), and in particular also the Jones-Wenzl relation, have to be
homogeneous with respect to any grading extending the original grading. Now assume that we are given $A$, $g_s, f_s$ as in (2) and that the Jones-Wenzl relations are homogeneous. We postpone the proof that the $H_3$ Zamolodchikov relation is homogeneous to Lemma \ref{lemma: homogeneity of H3} below. It is
straightforward to check that all the other relations \cite[(10.7a)-(10.7h), (10.14a)-(10.14d)]{EMTW} are homogeneous. (For example, if $\pa_s(v) \ne 0$ then $s(v) = v - \pa_s(v)
\alpha_s$, and the $W$-invariance of the grading implies that $\deg(v) = \deg(\alpha_s)$. Now one can confirm that \cite[(10.7d)]{EMTW} is homogeneous.) This shows that
\eqref{eq:multigrading} defines an $A$-grading on $\HC$ refining the original grading. \end{proof}

\begin{rem} The Zamolodchikov relations in types $A_3$, $B_3$, and $A_1 \times I_2(m)$ are homogeneous because the number of $2m_{s,t}$-valent vertices of each kind is the same on
both sides of the equality. The corresponding relation in type $H_3$ expresses the difference between two rex moves as a linear
combination of diagrams which include dots and trivalent vertices. This linear combination is not explicitly known. Proving its homogeneity will require a different argument, and we postpone it for pedagogical reasons. \end{rem}


\begin{cor} \label{cor:whengradingsagree} Two gradings on $\HC$ by an abelian group $A$ (which refine the original grading) agree if and only if they agree on the startdots of each color, and on $V$. \end{cor}

\begin{proof} By part (1) of Prop. \ref{prop: gradings extend if and only JW is homogeneous}, any $A$-grading is determined by $f_s$ and $g_s$ and the grading on $V$. Note that $f_s + g_s = \deg(\alpha_s)$, so $g_s$ is determined from $f_s$ and the grading on $V$. \end{proof}


Proposition \ref{prop: gradings extend if and only JW is homogeneous} reduces the problem of constructing gradings to verifying the homogeneity of the Jones-Wenzl relations. We now verify this homogeneity in the case where $f_s + g_s = f_t + g_t$. The general case is more difficult and will be discussed in \S\ref{sec-homog}.

\begin{lemma}\label{lemma: homogeneity of JW in the easy case}
    Assume we are given a $W$-invariant $A$-grading on $V$ and let $f_s,g_s,f_t,g_t \in A$ such that $f_s+g_s = f_t+g_t$. Then the Jones-Wenzl relation for $s,t$ is homogeneous (when generators have degrees as in \eqref{eq:multigrading}).
\end{lemma}
\begin{proof}
    Any term in the Jones-Wenzl idempotent is the deformation retract of a (two-colored) crossingless matching in some Temperley-Lieb algebra, see \cite[\S 5.3.2]{ECathedral}. Cups and caps deformation retract to pitchforks of various colors and we deduce
    \begin{subequations}
    \begin{equation} \deg(\pitchred) = f_t - g_s = \deg(\pitchblue) = f_s - g_t, \end{equation}
    \begin{equation} \deg(\pitchdownred) = g_t - f_s = \deg(\pitchdownblue) = g_s - f_t. \end{equation}
    \end{subequations}
    Thus the degree of (the deformation retract of) a cap is opposite that of a cup. Any crossingless matching with the same number of boundary points on top and bottom has the same number of cups as caps, so the overall degree of the deformation retract in $\HC$ is zero.
\end{proof}

We included the argument above because we find it edifying. One can also prove homogeneity of the Jones-Wenzl relation using an argument similar to Lemma \ref{lemma:
homogeneity of H3} below.

\begin{lemma} \label{lemma:secret inner} Suppose that $\HC$ is graded by $A$ and satisfies \eqref{eq:multigrading} for $f_s, g_s \in A$. Consider any diagram $D$ without
polynomials representing a morphism in $\HC$. For each $s \in S$, let $n_s$ be the number of copies of $s$ in the target object minus the number in the source object. Then there
exist integers $k_s$ for each $s \in S$ such that \begin{equation} \label{innerdegree} \deg(D) = \sum_{s \in S} n_s f_s + k_s(f_s + g_s). \end{equation} Moreover, the degree of $D$
in the original $\Z$-grading is $\sum_{s \in S} n_s + 2 k_s$. \end{lemma}

\begin{proof} This statement can be immediately verified for any identity map and for each of the generating morphisms in \eqref{eq:multigrading}. For example, if a generator has
one more $s$ in the target than the source, then there is a contribution of either $f_s$ or $-g_s$, and $-g_s = f_s - (f_s + g_s)$. Moreover, this statement is additive under composition or
tensor product of morphisms, whence the result follows. \end{proof}

Even though we do not know the $H_3$ Zamolodchikov relation explicitly, we know it must be homogeneous! The following lemma concludes the proof of Proposition \ref{prop: gradings extend if and only JW is homogeneous}.

\begin{lemma} \label{lemma: homogeneity of H3} Assume that there is a grading on $\HC$ valued in an abelian group $A$, which refines the original $\mathbb{Z}$-grading. Then any
Zamolodchikov relation of type $H_3$ is homogeneous (where generators have degrees as in \eqref{eq:multigrading}). \end{lemma}

\begin{proof} Let $m_{st} = 3$ and $m_{tu} = 5$ and $m_{su} = 2$ and consider the relation as on \cite[p11]{EWGr4sb}. Both sides of the relation are morphisms
	\[  (s,t,u,t,s,u,t,u,t,s,u,t,u,t,u) \to (u,t,u,t,u,s,t,u,t,u,s,t,u,t,s). \]
Both diagrams on the left-hand side are homogeneous of degree zero: for example, for each $6$-valent vertex of degree $g_s - g_t$ there is another of degree $g_t - g_s$. We need only show that the lower terms are also degree zero. They are linear combinations of diagrams which are themselves homogeneous, so we need only show that each diagram has degree zero. We note that any such diagram can be expressed without polynomials (but possibly with barbells), since the Zamolodchikov relation can be defined in a realization spanned by roots.

We pick a diagram $D$ without polynomials in this morphism space, and use the notation of the previous lemma. Note that $n_s = n_t = n_u = 0$. Choose $k_s \in \Z$ such that \eqref{innerdegree} holds. Since $\deg(D) = 0$ in the original $\Z$-grading, we must have $\sum k_s = 0$. Now \eqref{innerdegree} implies that
\begin{equation} \deg(D) = k_s(f_s + g_s) + k_t(f_t + g_t) + k_u(f_u + g_u). \end{equation}
If we can prove that $f_s + g_s = f_t + g_t = f_u + g_u$ in $A$, then $\deg(D)$ is zero as desired.

We claim that $f_s + g_s = f_t + g_t$. Otherwise, \eqref{eq:polygradingcriterion} implies that $[2]_{s,t} = [2]_{t,s} = 0$. But then $s$ and $t$ commute in their action on $V$, so
$(st)^2$ acts by the identity on $V$. Since $m_{st}$ is odd, this implies that $st$ acts by the identity on $V$. However, $st(\alpha_s) = - \alpha_s$, a contradiction. This same
argument applies for $t$ and $u$ since $m_{tu}$ is also odd. \end{proof}

\begin{proof}[Proof of Theorem \ref{thm:bigraded}] This follows immediately from  Proposition \ref{prop: gradings extend if and only JW is homogeneous} and Lemma \ref{lemma: homogeneity of JW in the easy case}. \end{proof}

\subsection{An annoying technicality related to the grading on polynomials} \label{subsec:polyssuck}

In the introduction we discussed that if $[2]_{s,t}$ or $[2]_{t,s}$ is nonzero, then $\alpha_s$ and $\alpha_t$ must have the same degree. Ignoring $\HC$ and focusing on $V$ for the
moment, they must have the same degree in any $W$-invariant grading on $V$ for which the simple roots are homogeneous. However, there are additional circumstances where $\alpha_s$ and $\alpha_t$ must have the same degree, even when $[2]_{s,t} = [2]_{t,s} = 0$, which all arise from Lemma \ref{lem:fdegree}. Let us restate that lemma in general, the proof being easily adapted from the introduction.

\begin{defn} Let $A$ be an abelian group. A \emph{good ($A$-valued) grading} on a realization $V$ is a $W$-invariant $A$-grading for which the simple roots are homogeneous. \end{defn}

\begin{lemma} \label{lem:fdegredux} Let $V$ be equipped with a good grading. If $s \in S$ and $f \in V$ are such that $s(f) \ne f$ and $f$ is homogeneous, then $\deg(f) = \deg(\alpha_s)$. \end{lemma}

\begin{ex} \label{ex:yuck} Consider a realization spanned by $\alpha_s$ and $\alpha_t$ and $f$, where $m_{st} = 2$, and $\pa_s(f) = \pa_t(f) = 1$. When $a, b \in \Z$, for reasons of
parity, any element of the form $g = f + a \alpha_s + b \alpha_t$ has the property that $\pa_s(g), \pa_t(g) \ne 0$. Thus if any such element $g$ is homogeneous, then $\deg(\alpha_s) =
\deg(g) = \deg(\alpha_t)$.

If the base ring of the realization is $\Z$ or the field $\mathbb{F}_2$, then in any good grading on $V$ there will be a homogeneous element of the form $g$ above.

If $2$ is invertible in the base ring, then $h := f - \frac{1}{2} \alpha_s - \frac{1}{2} \alpha_t$ is $W$-invariant. There is a good grading valued in $\Z^3$ where $\deg(\alpha_s) = (1,0,0) \ne \deg(\alpha_t) = (0,1,0)$, and $\deg(h) = (0,0,1)$. \end{ex}

As the previous example shows, the question of whether or not $\alpha_s$ and $\alpha_t$ must have the same degree in any good grading is quite subtle, and depends not just on the
Cartan matrix of the realization but on additional properties of the realization itself. The equivalence relation $\sim_V$, based only on the Cartan matrix, encodes a necessary but not
sufficient condition. In the next examples, we note that there are useful realizations where $\sim_V$ suffices.

\begin{ex} Suppose that $V$ has a basis consisting only of simple roots and $W$-invariant elements. Then there is a good grading where
$\deg(\alpha_s) = \deg(\alpha_t)$ if and only if $s \sim_V t$, and the $W$-invariant basis elements are homogeneous of arbitrary degrees. \end{ex}

\begin{ex} Suppose that the Cartan matrix of the realization is invertible over the base ring. Then we are in the situation of the previous example. \end{ex}

\begin{ex} \label{ex:universal}
 Fix a Cartan matrix associated to a Coxeter system $(W,S)$. The \emph{universal realization} associated to this Cartan matrix has a basis given by $\{\alpha_s, \varpi_s\}_{s \in S}$. The values of $\pa_s(\alpha_t)$ are determined by the Cartan matrix, and $\pa_s(\varpi_t) = \delta_{st}$ (the Kronecker delta). This is a realization satisfying Demazure
surjectivity (i.e. the map $\pa_s$ surjects from $V$ onto the base ring). Moreover, one can equip it with a good grading where $\deg(\alpha_s) = \deg(\varpi_s)$ for all $s$, and
$\deg(\alpha_s) = \deg(\alpha_t)$ if and only if $s \sim_V t$. \end{ex}

Again, we reiterate that, given a Cartan matrix, there exists a realization $V$ with that Cartan matrix having the following property: there exists a good grading such that
$\deg(\alpha_s) = \deg(\alpha_t)$ if and only if $s \sim_V t$. However, not every realization has this property. Similar statments hold for $\sim_C$. Ultimately, given a realization
$V$, we do not see an elegant criterion for when $\deg(\alpha_s) = \deg(\alpha_t)$ for \emph{every} good grading, beyond ad hoc extrapolation from Lemma \ref{lem:fdegredux}, nor do we
see great need for such a criterion. To make precise statements without worrying about this technicality, we adopt the approach that first one fixes a good grading on $V$, and then
extends it to a grading on $\HC$.

With these subtleties in mind, here is a generalization of the bigrading which (potentially) breaks color symmetry for different components.

\begin{thm}\label{thm: general grading by components}
Let $V$ be a realization of $(W,S)$. Equip $V$ with a good $W$-invariant grading, valued in an abelian group $\Gamma$, for which $\deg(\alpha_s) = \deg(\alpha_t)$ whenever $m_{st} \ne 2$. Let $\Lambda_{\CB}$ be as in Definition \ref{lambdadef}; equivalently, it is the free abelian group on generators $\{ f_s, g_s\}_{s \in S}$ modulo the relation $f_s + g_s = f_t + g_t$ whenever $m_{st} \ne 2$. Then the grading on $V \subset R$ extends to a grading on $\HC$ by the group $(\Gamma \times \Lambda_{\CB})/I$ where $I$ is the submodule generated by $(-\deg(\alpha_s),f_s + g_s)$ for all $s \in S$. The generating diagrams have degrees as in \eqref{eq:multigrading}.
\end{thm}

\begin{proof}
By Proposition \ref{prop: gradings extend if and only JW is homogeneous} it suffices to check that the Jones-Wenzl relation is homogeneous. For $m_{st} \neq 2$ this follows from Lemma \ref{lemma: homogeneity of JW in the easy case}. For $m_{st} = 2$ the Jones-Wenzl relation is clearly homogeneous since the corresponding Jones-Wenzl projector is trivial.
\end{proof}

\subsection{Inner gradings}

Consider a $\Z$-graded $\Bbbk$-linear category, and fix objects $X$ and $Y$. Let $Z = X(k)$ for some $k \in \Z$. Then $\Hom(X,Y)(-k) \cong \Hom(Z,Y)$ and $\Hom(Y,X)(k) \cong \Hom(Y,Z)$ as graded $\Bbbk$-modules. If we rename $Z$ as $X$, this is like equipping the category with a new
grading, where all the morphisms to and from $X$ are shifted by $\pm k$. One could rename any object in this fashion, or even rename the shift itself. The following is a small
generalization of this idea.

\begin{defn} Let $\CC$ be a $\Bbbk$-linear category which is graded by an abelian group $B$. An \emph{inner grading datum} is the data of:
	\begin{itemize}
		\item An abelian group $A'$ and a homomorphism $d \co B \to A'$,
		\item an element $d(X) \in A'$ for each object $X$ of $\CC$.
		\end{itemize}
If $\CC$ is monoidal then the inner grading datum is \emph{monoidal} if $d(\1) = 0$ and $d(X \otimes Y) = d(X) + d(Y)$. Given an inner grading datum, the corresponding \emph{inner grading} on $\CC$ is the $A'$-grading defined as follows: each morphism in $\Hom_{\CC}(X,Y)$ of $B$-degree $b$ is homogeneous in the new grading with $A'$-degree $d(Y) - d(X) + d(b)$.

An $A'$-grading is called \emph{inner} (relative to the original $B$-grading), if it corresponds to an inner grading datum. If $A \subset A'$ is a subgroup and all the degrees of
all morphisms live in the subgroup $A$, then we also abusively call the $A$-grading \emph{inner}. \end{defn}
	
\begin{example}\label{example: inner Frobenius grading} Let $A' \subset \Q \times \Q$ be the free $\Z$-module with basis $(1,0)$ and $(\frac{1}{2}, \frac{1}{2})$. Let us define a homomorphism $\Z \to A'$ sending $1
\mapsto (\frac{1}{2}, \frac{1}{2})$. For each $s \in S$, send the generator $B_s$ of $\HC$ to $d(B_s) = (\frac{1}{2}, -\frac{1}{2})$, and extend this function monoidally to a
function on all objects of $\HC$. This is an inner grading datum relative to the original $\Z$-grading on $\HC$, and the corresponding inner grading is the bigrading from Theorem
\ref{thm:bigraded}. Note that all degrees live inside the subgroup $A \subset A'$, where $A = \Z \times \Z$. \end{example}

\begin{remark} \label{rmk:maybeneedtoextendgroup} As in the previous example, it may not be possible to define the inner grading datum so that it is valued in a subgroup $A$ which contains all the degrees of
morphisms. \end{remark}

\begin{remark} Note that the collection of inner gradings is closed under specialization. \end{remark}

The following result is little more than a slightly more precise restatement of Lemma \ref{lemma:secret inner}. It proves that every grading is inner over a grading where $f_s = g_s$.

\begin{prop} \label{prop:inner} Assume that there is a grading on $\HC$ valued in an abelian group $A$, which refines the original $\mathbb{Z}$-grading. Let $\sim$ be the
equivalence relation defined by $s \sim t$ if and only if $\deg(\alpha_s) = \deg(\alpha_t)$ in $A$. Let $\Gamma \subset A$ be a subgroup such that $\deg(f) \in \Gamma$ for all
homogeneous elements $f \in V$. Then $\HC$ has a grading by the group $B := (\Gamma \times \Z^{S/\sim})/I$, where $I$ is the subgroup generated by $(-\deg(\alpha_s), 2_s)$ for all $s
\in S$. In this $B$-grading the degree of both univalent dots is $(0,1_s)$ and the degree of $f \in V$ is $(\deg(f),0)$. Finally, the original $A$-grading is inner over the
$B$-grading. \end{prop}

\begin{proof} Define $A'$ as the extension\footnote{The existence of $\For \colon A \to \Z$ eliminates worries about $2$-torsion in $A$.} of $A$ containing $\frac{f_s + g_s}{2}$ for all $s$. We define a homomorphism $d \co B \to A'$ sending $(0,1_s) \mapsto \frac{f_s + g_s}{2}$ for all $s \in S$, and sending $(\gamma,0)$ to $\gamma \in \Gamma \subset A$. This homomorphism clearly sends $I$ to zero, so it is well-defined. Set $d(B_s) = \frac{f_s - g_s}{2}$.  Take a diagram $D$ without polynomials, and compute the numbers $n_s$ and $k_s$ as in Lemma \ref{lemma:secret inner}. The degree of $D$ in the $\Z^{S/\sim}$ grading is $\sum_s (n_s + 2k_s) 1_s$. Note that $D$ represents a morphism $X \to Y$ for which $d(Y) - d(X) = \sum n_s (\frac{f_s - g_s}{2})$.  Thus the degree in $A'$ of the inner grading is
\begin{equation} \sum n_s (\frac{f_s - g_s}{2}) + (n_s + 2k_s) \frac{f_s + g_s}{2} = n_s f_s + k_s(f_s + g_s). \end{equation}
This degree lives in the subgroup $A \subset A'$, and by Lemma \ref{lemma:secret inner}, it agrees with $\deg(D)$ in the $A$-grading. The addition of polynomials to diagrams does not affect this calculation, since both groups treat $\Gamma$ the same. \end{proof}	

In particular, when $V$ is spanned by simple roots, or when the degrees of elements of $V$ are spanned by the degrees of simple roots, then we can set $\Gamma$ to be this span, and
hence $B \cong \Z^{S/\sim}$. Thus for a typical realization $V$ with a good grading, every grading on $\HC$ extending the grading on $V$ is inner over the grading by $\Z^{S/\sim}$.

%
%

\section{Autoequivalences of the Hecke category} \label{sec-autoequiv}

Proposition \ref{prop: gradings extend if and only JW is homogeneous} is also useful for studying autoequivalences.

\begin{defn}
Let $\CC$ be an additive $\Bbbk$-linear category which is graded by an abelian group $A$. Let $\chi : A \to \Bbbk^\times$ be a group homomorphism. Then $\CC$ admits an autoequivalence $\Theta_{\chi}$ which fixes all objects, and acts on homogeneous morphisms by the formula
\begin{equation} \Theta_{\chi}(f) = \chi(\deg(f)) \cdot f. \end{equation}
We call such an autoequivalence a \emph{degree-based rescaling}.
\end{defn}

It is straightforward to prove that degree-based rescalings are well-defined equivalences. Moreover, if $\CC$ is graded monoidal, then $\Theta_{\chi}$ is a monoidal equivalence,
as the monoidal composition is also compatible with the grading. If $\Theta$ is a degree-based rescaling then $\Theta$ acts diagonalizably on all Hom spaces. The converse is also true.

\begin{prop} Let $\CC$ be an additive $\Bbbk$-linear monoidal category with a $\Z$-grading. Let $\Theta$ be a monoidal, graded equivalence which fixes all objects and acts diagonalizably on all (homogeneous components of) morphism spaces. Then $\CC$ admits a grading by an abelian group $A$ which refines the original $\Z$-grading, and there is a homomorphism $\chi \co A \to \Bbbk^\times$ such that $\Theta = \Theta_{\chi}$. \end{prop}
	
\begin{proof} Let $X$ be the set of eigenvalues of $\Theta$, a subset of $\Bbbk^\times$ which contains $1$ (the eigenvalue of any identity map) and is closed under multiplication (one can multiply eigenvalues by taking tensor product of morphisms).  Let $A'$  be the subgroup of $\Bbbk^{\times}$ generated by $X$. Let $A = A' \times \Z$, where $\Z$ represents the original grading. Let $\chi \co A \to \Bbbk^\times$ be the projection to $A'$ followed by the inclusion map. Equip $\CC$ with a grading by $A$ by splitting each homogeneous piece of each hom space into $\Theta$-eigenspaces, where $A'$ records the eigenvalue and $\Z$ the original grading. It is clear that the grading is well-defined (i.e. composition and tensor product act appropriately on the grading) and that $\Theta = \Theta_{\chi}$. \end{proof}

\begin{thm} \label{thm:autoequiv} Let $\Theta$ be a monoidal graded autoequivalence of the (original, $\Z$-graded) category $\HC$ which fixes all objects and acts diagonalizably on $\End(\1)$. Then there is a grading on $\HC$ valued in an abelian group $A$, which refines the original $\mathbb{Z}$-grading, such that $\Theta = \Theta_{\chi}$ for some group homomorphism $\chi : A \to \Bbbk^\times$. In particular, two such autoequivalences agree if and only if they agree on the startdot and on $V$. \end{thm}

Before the proof, here is a useful and easy corollary.

\begin{cor} \label{cor:identitycriterion} Let $\Theta$ be a $\Z$-graded monoidal autoequivalence of $\HC$ which \begin{enumerate}
	\item fixes $B_s$ for each $s \in S$, 
	\item fixes polynomials, and
	\item fixes $\finaldotred$, the generator of $\Hom(B_s,\1)$ for each $s \in S$. \end{enumerate}
Then $\Theta$ is equal to the identity functor. \end{cor}

%

\begin{proof}[Proof of Theorem \ref{thm:autoequiv}]
Since $\startdotred$ lives in a one-dimensional graded morphism space, $\Theta$ sends $\startdotred$ to some scalar multiple of $\startdotred$. Write $\kappa_s$ for this scalar. Similarly, all non-polynomial generators of $\HC$ are sent  by $\Theta$ to scalar multiples of themselves. Write $\lambda_s$ for
the scalar of $\finaldotred$ under $\Theta$. By examining the unit and counit relations \cite[(8.5ab)]{EMTW} we see that the scalar on the split must be $\lambda_s^{-1}$, and the scalar on the merge must be $\kappa_s^{-1}$. In particular, by the barbell relation \cite[(10.7c)]{EMTW}, $\alpha_s$ is rescaled by $\kappa_s \lambda_s$. 

Because $\HC$ admits a set of generating morphisms which are rescaled by $\Theta$, we deduce that $\Theta$ rescales each diagram, and thus acts diagonalizably on each
finite-dimensional graded piece of each morphism space in $\HC$. (While we know that $\Theta$ acts diagonalizably on the span of the roots, it is a non-trivial assumption that
$\Theta$ acts diagonalizably on the rest of $V$.) Now we apply the previous proposition. Note that the group $A'$ of eigenvalues is the subgroup of $\Bbbk^\times$ generated by $\kappa_s$ and $\lambda_s$ and the eigenvalues of $\Theta$ on
$V$.

The final statement is clear from Corollary \ref{cor:whengradingsagree}. \end{proof}

\begin{rem} It may be helpful to elaborate on this proof, and think about classifying autoequivalences without using Proposition \ref{prop: gradings extend if and only JW is homogeneous}. Each relation of $\HC$, with the exception of \cite[(10.7d)]{EMTW}, is a linear combination of eigenvectors for $\Theta$. If they do not share the same eigenvalue, then one can derive a contradiction. For example, if $\mu$ is an eigenvalue of some diagram in the relation, then subtracting $\mu$ times the original relation from $\Theta$ applied to the relation, one gets a relation which does not involve that diagram. This relation does not hold in $\HC$. This is effectively the same argument that we used in the previous proof to argue that relations of $\HC$ were homogeneous for any grading refining the original grading.

Pick $s \ne t$. In order for the polynomial forcing relation \cite[(10.7d)]{EMTW} to be preserved by $\Theta$ when we set $f = \alpha_t$ we deduce that $\kappa_s \lambda_s = \kappa_t \lambda_t$ whenever $\pa_s(\alpha_t) \ne 0$. The reader should now see that the arguments involving $\Theta$ and the rescaling parameters like $\kappa_s$ are precisely analogous to the arguments involving homogeneity and degrees like $f_s$. \end{rem}

\begin{rem} The two-parameter family of autoequivalences coming from the bigrading of Theorem \ref{thm:bigraded} has been known to a handful of experts for some time. For example, a special case was used in \cite[Lemma 4.18]{MackaayTubbenhauer}. \end{rem}


By composing with a Dynkin diagram automorphism, this theorem also classifies autoequivalences which permute the generating objects $\{B_s\}$.


Corollary \ref{cor:identitycriterion} implies the following statement about functors out of the Hecke category.

\begin{cor} \label{cor:isomcriterion}  Let $F$ and $G$ be monoidal $\Z$-graded fully-faithful functors from $\HC$ to a graded additive monoidal category $\CC$, for which \begin{enumerate}
	\item we have chosen an isomorphism $\phi_1 \colon F(\one) \simto G(\one)$ which respects the monoidal identity structure,
	\item we have chosen isomorphisms $\phi_s \colon F(B_s) \simto G(B_s)$ for each $s \in S$, 
	\item $F$ and $G$ agree on polynomials, after intertwining with $\phi_1$, and
	\item $F$ and $G$ agree on $\finaldotred$ for each $s \in S$, after intertwining with $\phi_s$ and $\phi_1$. \end{enumerate}
Then $F \cong G$ via the unique monoidal natural transformation extending $\phi_1$ and $\{\phi_s\}$. \end{cor}

\begin{proof}
The idea of the proof is that any fully faithful functor can be inverted on its essential image, so, loosely speaking, we apply Corollary \ref{cor:identitycriterion} to $F^{-1} \circ G$. The more rigorous argument is as follows: For any object $X$ of $\HC$, let $\phi_X \co F(X) \to G(X)$ be the unique map which monoidally extends the maps $\phi_1$ and $\{\phi_s\}$. For example, if $X = B_s B_t$ then $\phi_X = \phi_s \ot \phi_t$. Any object is a tensor power of $\one$ and various $B_s$, so this map $\phi_X$ makes sense. Let $\Theta_{X,Y} \co \Hom_{\HC}(X,Y) \to \Hom_{\HC}(X,Y)$ be the unique linear map making the following diagram commute.
\begin{equation}
	\begin{tikzcd}[column sep = huge]
	\Hom_{\HC}(X,Y) \arrow[d,"F"] \arrow[r,"\Theta_{X,Y}"] & \Hom_{\HC}(X,Y) \arrow[d,"G"] \\
	\Hom_{\CC}(F(X),F(Y)) \arrow[r,"\phi_Y \circ (-) \circ \phi_X^{-1}"] & \Hom_{\CC}(G(X),G(Y))
	\end{tikzcd}
\end{equation}
By assumption, $\Theta_{\one, \one}$ is the identity map on polynomials, and $\Theta_{B_s,\one}$ fixes the dot. Observe that $\Theta_{X,Y}$ is the action on morphism spaces of some autoequivalence $\Theta$ which fixes objects. Then Corollary \ref{cor:identitycriterion} finishes the proof.
\end{proof}

Both corollaries are tools which allow one to prove results about autoequivalences or functors with the minimum of computational effort. For example, using Corollary
\ref{cor:identitycriterion} we can avoid computing how $\Theta$ acts on most of the generators of the Hecke category, including the complicated $2m$-valent vertices.

\begin{rem} Corollary \ref{cor:isomcriterion} was applied in \cite{EHDrinfeld} to give quick proofs of folklore results about conjugation by braids. Let $K^b(\HC)$ denote the category of bounded chain complexes in $\HC$, with morphisms considered modulo homotopy. For example, it is a folklore theorem for $W =
S_n$ that the (Rouquier complex of the) full twist braid $FT$ is in the Drinfeld center of the Hecke category, i.e. conjugation by the full twist is isomorphic to the identity
functor of $K^b(\HC)$. Using Corollary \ref{cor:isomcriterion}, one can efficiently prove that the two functors $\HC \to K^b(\HC)$ below, \begin{equation} X \mapsto X \qquad \text{
and } \qquad X \mapsto FT \otimes X \otimes FT^{-1}, \end{equation} are isomorphic.  It is less well-known that conjugation by the half twist braid is equivalent to the composition of a Dynkin diagram automorphism, and a degree-based rescaling $\Theta_{\chi}$. Here we use the bigrading of Theorem \ref{thm:bigraded}, and $\chi(a,b) = (-1)^b$. \end{rem}

\begin{rem} The previous remark treats isomorphisms between functors $\HC \to K^b(\HC)$. These functors extend to functors $K^b(\HC) \to K^b(\HC)$, though it does not automatically follow that the natural isomorphism of functors will extend to the homotopy category. We prove this extension result in \cite{EHDrinfeld} as well. \end{rem}

\section{The Frobenius automorphism} \label{sec:frobenius}

This section is meant to justify why the bigrading has anything to do with the Frobenius operator in geometry. We give two separate explanations, one for the reader familiar with geometry, and another one purely algebraic.

\subsection{Purity in geometry}

We recall the geometric incarnation of the Hecke category. Let $G$ be a reductive group with Borel subgroup $B$, flag variety $\mathcal{B}$ and Weyl group $G$. The derived equivariant category $D^b_B(\mathcal{B})$ of constructible sheaves can be equipped with a convolution product $\star$. For any $s \in S$ denote by $\ICs_s = \textbf{1}_{Y_s} [1]$ the simple perverse sheaf supported on the Schubert variety $Y_s = \overline{BsB}/B$.  Objects of the form $IC_{s_1} \star... \star \ICs_{s_k}$ are called Bott-Samelson sheaves. For objects $X, Y \in D^b_B(\mathcal{B})$ we write $\Hom^i(X,Y) := \Hom(X,Y[i])$. Let $\HC_{geo}$ be the monoidal category with objects all Bott-Samelson sheaves and with (graded) morphism spaces
\begin{equation*}
\Hom_{\HC_{geo}}(X,Y) = \bigoplus_{i \in \mathbb{Z}} \Hom^i_{D^b_B(\mathcal{B})}(X, Y).
\end{equation*}
Then there is an equivalence of (graded) monoidal categories
\begin{align*}
\HC \overset{\sim}{\rightarrow} \HC_{geo}, \qquad B_s \mapsto \ICs_s.
\end{align*}

Assume now that $G$ is defined in positive characteristic. Then $\ICs_s$, and thus also any Bott-Samelson sheaf, can be equipped with a canonical $\mathbb{F}_q$-structure (here we follow the convention that the canonical $\mathbb{F}_q$-structure of $\ICs_s$ is pure of weight $1$). This equips the morphism space between any two Bott-Samelson sheaves with a canonical Frobenius action which can be computed explicitly, and the following lemma states that this action is pure.

\begin{lemma}\label{lemma: Frobenius purity}
	The action of Frobenius on $\Hom^i (\ICs_{s_1} \star ... \star\ICs_{s_k} , \ICs_{s_1'} \star... \star \ICs_{s_l'})$ given by multiplication with $q^{\frac{i+l-k}{2}}$ and the Hom-space vanishes unless $i+l-k$ is even.
\end{lemma}
\begin{proof}
We briefly sketch the argument of this standard result for the convenience of the reader. We may first use that $\Hom( X \star \ICs_s , Y) \cong \Hom(X, Y \star \ICs_s)$ (see \cite[Proposition 7.6.10]{AcharPerv}) to reduce to the case where $l = 0$. In this case, we we want to compute the Frobenius action on
\begin{equation*}
	\Hom^i_{D^b_B(\mathcal{B})} (\ICs_{s_1} \star...  \star\ICs_{s_k}, \textbf{1}_{eB/B})  = \Hom^i_{D^b_B(pt)}(( \ICs_{s_1} \star...  \star\ICs_{s_k})_{eB/B}, \textbf{1}_{eB/B}).
\end{equation*} 
It is known that $H^i(\ICs_{s_1} \star...  \star\ICs_{s_k})_{eB/B}$ is a direct sum of objects of the form $\textbf{1}_{eB/B} ( -\frac{k+ i}{2})$ for $i+k$ even and vanishes for $i+k$ odd (see \cite[Proposition 7.3.6, Lemma 7.4.7]{AcharPerv}). A standard long exact sequence argument then proves that Frobenius acts on $\Hom^i_{D^b_B(pt)}(( \ICs_{s_1} \star...  \star\ICs_{s_k})_{eB/B}, \textbf{1}_{eB/B})$ by multiplication with $q^{\frac{i-k}{2}}$ and that this space vanishes for $i-k$ odd. 
\end{proof}

We can refine the $\mathbb{Z}$-grading on $\HC_{geo}$ to a $\mathbb{Z}^2$-grading where the $(i,j)$ component is the $q^j$-eigenspace in $\Hom^i$. This yields
\begin{equation} \deg \left( \startdotred \right) = (1,1), \quad \deg \left( \finaldotred \right) = (1,0), \quad \deg \left( \splitred \right) = (-1,0), \quad 
	\deg \left( \mergered \right) = (-1,-1), \end{equation}
\begin{equation*} \deg \left( \sbotttop \right) = (0,0). \end{equation*}

We can also rewrite this grading along the isomorphism
\begin{equation} \label{eq:regradingmaybe}
\mathbb{Z}^2 \rightarrow \mathbb{Z}^2, (n,m) \mapsto n\cdot (0,1) + m \cdot ( 1, -1 )
\end{equation}
to obtain the bigrading from Theorem \ref{thm:bigraded}. We have already seen that this grading is inner to the standard $\mathbb{Z}$-grading. However, this also follows immediately from the purity result in Lemma \ref{lemma: Frobenius purity}. Thus, inner gradings arise naturally from purity in geometry.


\subsection{An algebraic Frobenius}

We thank Olivier Dudas for explaining this perspective to us, see also Remark \ref{rmk:thanksOlivier}. We are brief in our justifications, but the calculations are easy.

Let $q$ be an invertible element of the base ring $\Bbbk$, possibly a formal parameter. Let $F \co V \to V$ be multiplication by $q$. Extend $F$ to an automorphism of $R$. Then,
for the ordinary grading on $R$, $F$ just rescales each homogeneous element by a power of $q$ according to its degree. One can view $F \co R \to R$ as the action of Frobenius on the
equivariant cohomology of a point.

Given any $R$-bimodule $M$, we can define a new $R$-bimodule $F^{*}(M)$ with the same underlying $\Bbbk$-module, as follows. If $m \in M$, then we also let $m$ denote the corresponding element in $F^*(M)$, but we use $\cdot_F$ for the action on $F^*(M)$.  Then for $f, g \in R$ we have 
\begin{equation} f \cdot_F m' \cdot_F g := F(f) \cdot m \cdot F(g). \end{equation}
The operation $F^*$ is functorial, sending a bimodule morphism to the bimodule morphism with the same underlying linear map. 
The functor $F^*$ is monoidal. Letting $\1$ denote the regular $R$-bimodule $R$, we have $F^*(\1) \cong \1$. Since $F(R^s) = R^s$ we also have $F^*(B_s) \cong B_s$. Thus $F^*$ preserves all Bott-Samelson bimodules up to isomorphism. Henceforth we fix a monoidal family of such isomorphisms $\phi$ so that they preserve the one-tensor. For example, 
\begin{equation} \phi_s \co B_s \to F^*(B_s), \qquad \phi_s(f \ot g) = F(f) \ot F(g). \end{equation}

\begin{lemma} There are commuting diagrams
\begin{equation} \label{dotsquares} \begin{tikzcd}
B_s \arrow{r}{\finaldotred} \arrow[swap]{d}{\phi_s} & R \arrow{d}{\phi_{\one}} \\
F^*(B_s) \arrow{r}{F^*(\finaldotred)} & F^*(R)
\end{tikzcd} \qquad
\begin{tikzcd}
R \arrow{r}{\startdotred} \arrow[swap]{d}{\phi_{\one}} & B_s \arrow{d}{\phi_s} \\
F^*(R) \arrow{r}{qF^*(\startdotred)} & F^*(B_s)
\end{tikzcd}.
\end{equation}
\end{lemma}

For example, both paths around the first square send $1 \ot 1 \in B_s$ to $1 \in F^*(R)$, while both paths around the second square send $1 \in R$ to $\frac{q}{2}(\alpha_s \ot 1 + 1 \ot \alpha_s) \in F^*(B_s)$. Note that \eqref{dotsquares} does not commute without the extra factor of $q$ in one square, and one should not expect it to; the functor $F^*$ is not isomorphic to the identity functor, and $\phi$ is not a natural transformation.

Using the maps $\phi_M$ to identify $M$ with $F^*(M)$, we can produce a monoidal functor $G$ isomorphic to $F^*$ which sends $B_s$ to $B_s$ on the nose (not just up to
isomorphism). If we want the square
\begin{equation} \begin{tikzcd}
R \arrow{r}{G(\startdotred)} \arrow[swap]{d}{\phi_{\one}} & B_s \arrow{d}{\phi_s} \\
F^*(R) \arrow{r}{F^*(\startdotred)} & F^*(B_s)
\end{tikzcd}
\end{equation}
to commute, so that $\phi$ is an isomorphism from $G$ to $F^*$, then we need $G(\startdotred) = q^{-1} \startdotred$.

\begin{thm} The autoequivalence $G$ comes from an autoequivalence of the diagrammatic Hecke category, satisfying
\begin{equation} G(\startdotred) = q^{-1} \startdotred, \quad G(\finaldotred) = \finaldotred, \quad G(v) = q^{-1} v \text{ for all } v \in V. \end{equation}
In particular, $G$ is a degree-based rescaling associated to the bigrading in Theorem \ref{thm:bigraded}, via the homomorphism $\chi : \Z^2 \to \Bbbk^\times$ sending $(1,0) \mapsto q^{-1}$ and $(0,1) \mapsto 1$.
\end{thm}

\begin{rem} \label{rmk:thanksOlivier} The Hecke category has an object-fixing autoequivalence constructed using conjugation by the half twist Rouquier complex, composed with the
associated Dynkin diagram automorphism. Matt Hogancamp and the second author computed this autoequivalence, and obtained $G$ specialized at $q=-1$ (one can now find the proof in
\cite{EHDrinfeld}, based on the results of this paper). In conversations with Olivier Dudas at WARTHOG 2018, he explained to the second author that conjugation by the half twist is
supposed to agree with the Frobenius autoequivalence at $q=-1$. For this reason, the second author came up with the autoequivalence $G$, and conjectured it to be the
Soergel-ification of the Frobenius autoequivalence, in an earlier version \cite{EBigraded} of this paper.  \end{rem}

\begin{rem} The Soergel functor from $D^b_B(\mathcal{B})$ to $R$-bimodules intertwined the homological shift with the usual grading shift. This chapter suggests that one should define a mixed version of the Soergel functor which intertwines the algebraic and geometric bigradings. We recommend using the regrading \eqref{eq:regradingmaybe} when passing from geometry to algebra. \end{rem}

\section{Homogeneity of the Jones-Wenzl relation} \label{sec-homog}

Suppose we have a realization $V$ with a good grading. When $\langle \alpha_s, \alpha_t^{\vee} \rangle = \langle \alpha_t, \alpha_s^{\vee} \rangle = 0$ it can happen that $\alpha_s$ and $\alpha_t$ have distinct degrees, as in the following prototypical example.

\begin{example}\label{example: dihedral bigraded realization}
    Let $W = \langle s,t  \mid s^2 = t^2 = (st)^{m_{st}}  = 1\rangle$ be a dihedral group with $m_{st} = 2p^k$ where $p \neq 2$ is a prime. Let $\Bbbk$ be a ring with $p = 0$. Then $V = \Bbbk \alpha_s \oplus \Bbbk \alpha_t$ is a realization of $W$ where $\langle \alpha_s, \alpha_t^{\vee} \rangle = \langle \alpha_t, \alpha_s^{\vee} \rangle = 0$. We can equip $V$ with a $W$-invariant bigrading by setting $\deg(\alpha_s) = (2,0)$ and $\deg(\alpha_t) = (0,2)$.
\end{example}

\subsection{Proof of homogeneity}

We work in the setting of Proposition \ref{prop: gradings extend if and only JW is homogeneous}, i.e. we assume we are given an abelian group $A$ and a $W$-invariant $A$-grading on $V$ together with elements $g_s, f_s \in A$ such that \eqref{eq:multigrading} and \eqref{eq:polygradingcriterion} hold. Pick a pair $s \ne t \in S$. We have already seen that the Jones-Wenzl relation is homogeneous whenever $f_s+g_s = f_t +g_t$. Let us assume that $f_s + g_s \neq f_t + g_t$. Then we necessarily have $\langle \alpha_s, \alpha_t^{\vee} \rangle = \langle \alpha_t, \alpha_s^{\vee} \rangle = 0$, so the Jones-Wenzl idempotent comes from the Temperley-Lieb algebra $ \TL_{m_{st}-1}(0)$ on $m_{st}-1$ strands over $\Bbbk$ specialized at $[2]= 0$. Since $V$ is a realization of $(W,S)$ in the sense of \cite[Theorem B]{HaziRotatable}, the Jones-Wenzl projector $\JW \in \TL_{m_{st}-1}(0)$ exists and is rotatable. This puts the following restriction on $m_{st}$.

\begin{lemma}\label{lemma: m is 2 p^k}
Suppose $2 < m_{st} < \infty$ and $\Bbbk$ is a field. The Jones-Wenzl projector exists and is rotatable in $TL_{m_{st}-1}(0)$ if and only if $\Bbbk$ has positive characteristic $p$ and $m_{st} = 2p^k$ for some $k>0$.
\end{lemma}

In the following proof ${a \brack b}$ represents a quantum binomial coefficient, living in the polynomial ring $\Z[\delta]$ where $\delta = [2]$. Meanwhile ${a \choose b}$
represents an ordinary binomial coefficient, living in $\Z$. Both can be viewed as elements of $\Bbbk$ after specializing $[2]$ to zero. Note that ratios like $\frac{[6]}{[2]}$ are
also elements of $\Z[\delta]$ and have well-defined specializations, even though the numerator and denominator both specialize to zero. However, ratios like $\frac{[10]}{[6]}$ need
not make sense in $\Bbbk$, since they live in the fraction field of $\Z[\delta]$ but not in $\Z[\delta]$ itself.

\begin{proof}
Consider the map $\Z[\delta] \to \Z$ sending $\delta \mapsto 0$. For any $a, b \ge 0$ we claim that	
\begin{subequations} \label{quantumbinomclaim}
    \begin{align}
        {2a \brack 2b } &\mapsto  {a \choose b},\\
        {2a+1 \brack 2b } &\mapsto (-1)^{b} {a \choose b},\\
        {2a \brack 2b+1 } &\mapsto 0,\\
        {2a+1 \brack 2b+1 } &\mapsto (-1)^{a+b} {a \choose b},
    \end{align}
\end{subequations}
We prove this claim after using it to finish the proof.
	
By \cite{HaziRotatable} existence and rotatability is equivalent to the quantum binomial coefficient ${m_{st} \brack k }$ vanishing in $\Bbbk$ for all $1 \le k \le m_{st}-1$. If $m_{st} = 2d+1$ is odd then ${m_{st} \brack m_{st}-1} = (-1)^{d}{d \choose d} \ne 0$. So ${m_{st} \brack k }$ vanishes for all $1 \le k \le m_{st}-1$ if and only if $m_{st} = 2d$ and ${d \choose b}  = 0$ for all $1 \le b < d$. This is equivalent to $d$ being a power of $p$ in positive characteristic $p$ (or $d = 1$ in characteristic $0$).

It remains to prove \eqref{quantumbinomclaim}. This is a computation in the ring $\Z \cong \Z[\delta]/(\delta)$. However, there is a commuting square of rings
\begin{equation}
\begin{tikzcd}
{\Z[\delta]} \arrow[d] \arrow[r, hook] & {\Z[q,q^{-1}]} \arrow[d] \\
\Z \arrow[r, hook]                     & {\Z[\sqrt{-1}]}         
\end{tikzcd}
\end{equation}
where one thinks of $\Z[\sqrt{-1}]$ as $\Z[q]/(q^2+1)$. As the map $\Z \to \Z[\sqrt{-1}]$ is injective, we can prove the result in $\Z[\sqrt{-1}]$ instead. As a consequence, we can work before specializing in $\Z[q,q^{-1}]$ where $[2]=q+q^{-1}$.

Applying the quantum binomial theorem (c.f. \cite[1.3.1 (c)]{LuszQuantum}) with $q = \sqrt{-1}$, we get
\begin{equation*}
    \prod_{j=0}^{u-1} (1+\sqrt{-1}^{2j}z) = \sum_{t=0}^u \sqrt{-1}^{t(u-1)} {u \brack t} z^t
\end{equation*}
Let us first take $u = 2a$. Then the left hand side simplifies to
\begin{equation*}
    \prod_{j=0}^{2a-1} (1+(-1)^j z) = (1-z^2)^a.
\end{equation*}
Comparing the coefficient of $z^{2b+1}$ immediately yields $ {2a \brack 2b+1} = 0$. Comparing the coefficient of $z^{2b}$ we get $ (-1)^b{a \choose b}  =\sqrt{-1}^{2b(2a-1)} {2a \brack 2b} = (-1)^b {2a \brack 2b}$ and thus ${2a \brack 2b} = {a \choose b}$. For $u = 2a+1$ the left hand side becomes
\begin{equation*}
    \prod_{j=0}^{2a} (1+(-1)^j z) = (1+z)(1-z^2)^a
\end{equation*}
Comparing the coefficient of $z^{2b}$ yields $(-1)^b {a \choose b} = \sqrt{-1}^{2b \cdot 2a} {2a+1 \brack 2b} = {2a+1 \brack 2b}$. Comparing the coefficient of $z^{2b+1}$ yields $ (-1)^b {a \choose b} = \sqrt{-1}^{(2b+1) \cdot 2a} {2a+1 \brack 2b+1} = (-1)^a {2a+1 \brack 2b+1}$ and thus ${2a+1 \brack 2b+1} = (-1)^{a+b} {a \choose b}$.
\end{proof}

\begin{example} To put \eqref{quantumbinomclaim} in context it helps to learn some basic facts about quantum numbers after specializing $[2]=0$. The classic formula $[2][n] = [n+1]+[n-1]$ implies by an easy induction that
\begin{equation} \label{eq:quantumwhen2iszero} [2a] = 0, \qquad [2a+1] = (-1)^{a} \end{equation} for all $a \ge 0$.
Using \eqref{eq:quantumwhen2iszero} one can compute, for example, that
\begin{equation} \frac{[6]}{[2]} = [5] - [3] + [1] = 1 - (-1) + 1 = 3. \end{equation}
More generally, since
\begin{equation} \frac{[2a]}{[2]} = [2a-1] - [2a-3] + [2a-5] - \ldots \pm [1], \end{equation}
we deduce that
\begin{equation} \label{eq:2aover2} \frac{[2a]}{[2]} = (-1)^{a-1} a. \end{equation}
More generally, if $b$ divides $a$, then we have
\begin{equation} \frac{[2a]}{[2b]} = [2a-2b+1] - [2a-2b-1] + [2a-6b+1] - [2a-6b-1] + \ldots, \end{equation}
an expression with $\frac{a}{b}$ total terms that either ends in $+[2b+1] - [2b-1]$ or ends in $+[1]$, depending on parity. Regardless, we deduce that
\begin{equation} \frac{[2a]}{[2b]} = (-1)^{a-b} \frac{a}{b}. \end{equation}

So for example one has
\begin{equation} {8 \brack 4} = \frac{[8][7][6][5]}{[4][3][2][1]} = \frac{[8]}{[4]} \frac{[6]}{[2]} \frac{(-1)(1)}{(-1)(1)} = \frac{4 \cdot 3}{2 \cdot 1} = {4 \choose 2}. \end{equation}
\end{example}

Let us abbreviate $m_{st}$ to $m$. Suppose that $\Bbbk$ has characteristic $p$ and $m = 2p^k$ for some $k \ge 1$. By a \emph{TL diagram} we mean a crossingless matching on the planar strip with $m-1$ endpoints on top and bottom. For any $X \in \TL_{m-1}(0)$, we denote by
\begin{equation*}
\HeckX \in \End_{\HC_{\uni}}( sts... t)
\end{equation*}
the corresponding morphism in the universal two-color diagrammatic Hecke category $\HC_{\uni}$ of $s$ and $t$ obtained by deformation retract, see \cite[\S 5.3.2]{ECathedral}. Note that by Proposition \ref{prop: gradings extend if and only JW is homogeneous} that the $A$-grading on $V$ extends to $\HC_{\uni}$ (since there is no Jones-Wenzl relation). If $\HeckX$ is homogeneous, we say that $X$ is homogeneous and define the degree of $X$ to be
\begin{equation*}
\deg(X) := \deg \left( \HeckX \right) \in A.
\end{equation*}
Note that TL diagrams are always homogeneous in the original grading but may have non-zero $A$-degree. For example, we have
\begin{equation*}
\deg \left( \vcenter{\hbox{\begin{tikzpicture}[baseline={(current bounding box.center)}, scale=0.1]
    \draw [dashed] (-4,3) -- (4,3);
    \draw [dashed] (-4,-3) -- (4,-3);
    \draw (2.5,-3) -- (2.5,3);
    \draw [black] (1,-3) arc (0:180:2);
    \draw [black] (1,3) arc (0:-180:2);
\end{tikzpicture}}}\right) =
\deg \left( \vcenter{\hbox{\begin{tikzpicture}[baseline={(current bounding box.center)}, scale=0.07]
    \draw [dashed] (-7,5) -- (7,5);
    \draw [dashed] (-7,-5) -- (7,-5);
    \draw[red, line width=0.05cm] (-1,-0.5) -- (-5,-5);
    \draw[red, line width=0.05cm] (-1,-0.5) -- (3,-5);
    \draw[red, line width=0.05cm] (-1,0.5) -- (-5,5);
    \draw[red, line width=0.05cm] (-1,0.5) -- (3,5);
    \draw[red, line width=0.05cm] (-1,0.5) -- (-1,-0.5) ;
    \draw[blue, line width=0.05cm] (-1,-5) -- (-1,-3);
    \fill[blue] (-1,-3) circle (0.75cm);
    \draw[blue, line width=0.05cm] (-1,5) -- (-1,3);
    \fill[blue] (-1,3) circle (0.75cm);
    \draw[blue, line width=0.05cm] (5,-5) -- (5,5);
\end{tikzpicture}}}\right)
=f_t+g_t-f_s-g_s = \deg(\alpha_t) - \deg(\alpha_s)
\end{equation*}
which may be non-zero (see Example \ref{example: dihedral bigraded realization}).

To determine the degree of a TL diagram, one may proceed as follows. Any cap in a diagram whose left end is the $i$-th point on bottom
\begin{equation*}
\vcenter{\hbox{\begin{tikzpicture}[baseline={(current bounding box.center)}, scale=0.2]
    \draw [dashed] (-5,-5) -- (5,-5);
    \draw [black] (4,-5) arc (0:180:4);
    \draw [black] (-2,-5) rectangle (2,-3);
    \node at (0,-4.3) {\textbf{*}};
    \node[font=\tiny] at (-4,-5.5) {$i$};
\end{tikzpicture}}}
\end{equation*}
locally deforms in $\HC_{\uni}$ to
\begin{equation*}
    \vcenter{\hbox{\begin{tikzpicture}[baseline={(current bounding box.center)}, scale=0.2]
    \draw [dashed] (-5,-6) -- (5,-6);
    \draw [red, line width=0.4mm] (0,-3) -- (0,-2);
    \draw [red, line width=0.4mm] (-1.5,-6) -- (-1.5,-5);
    \draw [blue, line width=0.4mm] (-1,-6) -- (-1,-5);
    \draw [red, line width=0.4mm] (-0.5,-6) -- (-0.5,-5);
    \draw [red, line width=0.4mm] (1.5,-6) -- (1.5,-5);
    \draw [blue, line width=0.4mm] (5,-6) -- (0,0);
    \draw [blue, line width=0.4mm] (-5,-6) -- (0,0);
    \draw [blue, line width=0.4mm] (0,0) -- (0,2);
    \draw [black] (-2,-5) rectangle (2,-3);
    \fill[red] (0,-2) circle (0.3cm);
    \node at (0,-4.3) {\textbf{*}};
    \foreach \x in {0,0.5,1} {
        \node[circle, fill, inner sep=0.01cm] at (\x,-5.5) {};
    }
\end{tikzpicture}}}
\end{equation*}
if $i$ is even and to
\begin{equation*}
\vcenter{\hbox{\begin{tikzpicture}[baseline={(current bounding box.center)}, scale=0.2]
    \draw [dashed] (-5,-6) -- (5,-6);
    \draw [blue, line width=0.4mm] (0,-3) -- (0,-2);
    \draw [blue, line width=0.4mm] (-1.5,-6) -- (-1.5,-5);
    \draw [red, line width=0.4mm] (-1,-6) -- (-1,-5);
    \draw [blue, line width=0.4mm] (-0.5,-6) -- (-0.5,-5);
    \draw [blue, line width=0.4mm] (1.5,-6) -- (1.5,-5);
    \draw [red, line width=0.4mm] (5,-6) -- (0,0);
    \draw [red, line width=0.4mm] (-5,-6) -- (0,0);
    \draw [red, line width=0.4mm] (0,0) -- (0,2);
    \draw [black] (-2,-5) rectangle (2,-3);
     \fill[blue] (0,-2) circle (0.3cm);
    \node at (0,-4.3) {\textbf{*}};
    \foreach \x in {0,0.5,1} {
    \node[circle, fill, inner sep=0.01cm] at (\x,-5.5) {};
    }
\end{tikzpicture}}}
\end{equation*}
if $i$ is odd. Hence, the cap contributes $g_s-f_t$ or $g_t-f_s$ to the total degree of the TL diagram depending on the parity of $i$. Similarly, a cup whose left end point is the i-th point on top contributes either $f_s-g_t$ or $f_t-g_s$.
\begin{definition}
A cup/cap in a TL diagram is called even (resp. odd) if the index of its left endpoint is even (resp. odd). Here we index both the top and the bottom rows by $\{1,...,m-1\}$ reading from left to right.
\end{definition}
Let $D$ be a TL diagram. The local degree description for cups/caps from above gives rise to the following degree formula:
\begin{equation}\label{eq: degree formula}
    \deg(D) = (\# \{\text{even caps}\}- \# \{ \text{odd cups} \}) \cdot (g_s-f_t) + (\# \{\text{odd caps} \} -  \# \{ \text{even cups}\}) \cdot (g_t-f_s). 
\end{equation}
Note that in any TL diagram, the number of cups equals the number of caps and thus
\begin{equation*}
   \# \{\text{even caps}\}- \# \{ \text{odd cups} \} = -(\# \{\text{odd caps} \} -  \# \{ \text{even cups}\} ).
\end{equation*}
Hence, the degree of a TL diagram can also be written as
\begin{equation}\label{eq: degree formula2}
\begin{aligned}
    \deg(D) & = (\# \{\text{even caps}\}- \# \{ \text{odd cups} \}) \cdot ((g_s+f_s)-(g_t +f_t)). \\
\end{aligned}
\end{equation}

\begin{lemma}\label{lem: rotation introduces minus on degree}
Let $D$ be a TL diagram and let $D'$ be the diagram obtained by (counterclockwise) rotation by one strand. Then $\deg(D) = - \deg(D')$.
\end{lemma}
\begin{proof}
By \eqref{eq: degree formula2} it suffices to shows that
\begin{equation}\label{eq: degree formula for rotation}
    \# \{\text{even caps in }D\}- \# \{ \text{odd cups in }D \} =  \# \{\text{odd caps in } D' \} -  \# \{ \text{even cups in } D'\}.
\end{equation}
Rotation swaps the parity of cups (resp. caps) and preserves through strands, except possibly for the strand $A$ containing the first point in the top row and the strand $B$ containing the last point in the bottom row of $D$. Thus, it suffices to show that the combined contribution of $A$ and $B$ to the left hand side of \eqref{eq: degree formula for rotation} equals the contribution of the corresponding rotated strands $A'$ and $B'$ on the right hand side of \eqref{eq: degree formula for rotation}. If $A$ and $B$ are the same strand, then $A = B$  is a through strand which gets transformed by rotation into another through strand $A' =B'$. Thus, the contribution is $0$ on both sides. If $A \neq B$ there are four possibilities:

{\bf Case 1}: $A$ and $B$ are both through strands. Then their contribution is $0$ on the left hand side. $A'$ is an odd cap on the bottom left and $B'$ is an even cup on the top right (using that $m-1$ is odd). Their combined contribution in the right hand side is still $0$.

{\bf Case 2}: $A$ is a cup and $B$ is a cap. Then $A$ is odd and $B$ is even which have a combined contribution of $0$. $A'$ and $B'$ are both through strands which do not contribute to the right hand side.

{\bf Case 3}: $A$ is a cup and $B$ is a through strand. Then $A$ and $B$ contribute $-1$ to the left hand side since $A$ is odd. $A'$ is a through strand and $B'$ cup which also contribute $-1$ to the right hand side since $B'$ is even.

{\bf Case 4}: $A$ is a through strand and $B$ is a cap. Then $A$ and $B$ contribute $1$ to the left hand side since $B$ is even. $A'$ is a cap and $B'$ is a through strand which also contribute $1$ to the right hand side since $A'$ is odd.

We have now verified \eqref{eq: degree formula for rotation} in all cases which completes the proof.
\end{proof}

For any Temperley-Lieb diagram $D$, write $\coeffJW D \in \Bbbk$ for the coefficient of $D$ in the Jones-Wenzl projector $\JW$. When we draw a diagram we will sometimes use labeled strands for multiple parallel strands and nested cups/caps. For example,
\begin{equation*}
    \vcenter{\hbox{\begin{tikzpicture}[baseline={(current bounding box.center)}, scale=0.1]
    \draw [dashed] (-11,5) -- (1,5);
    \draw [dashed] (-11,-5) -- (1,-5);
    \draw [dashed] (-11,5) --(-11,-5);
    \draw [dashed] (1,5) -- (1,-5);
    \draw (-5,-5) -- (-5,5);
    \draw [black] (-6,5) arc (0:-180:2);
    \draw [black] (0,-5) arc (0:180:2);
    \draw [black] (-6,-5) arc (0:180:2);
    \node [font = \tiny]at (-8,-2) {3};
    \node [font = \tiny]at (-4,0) {2};
    \node [font = \tiny]at (-8,2) {2};
    \draw [black] (0,5) arc (0:-180:2);
    \node [font = \tiny]at (-1,2) {2};
    \end{tikzpicture}}}
    \text{ } := \text{ }
    \vcenter{\hbox{\begin{tikzpicture}[baseline={(current bounding box.center)}, scale=0.1]
    \draw [dashed] (-11,5) -- (6,5);
    \draw [dashed] (-11,-5) -- (6,-5);
    \draw [dashed] (-11,5) --(-11,-5);
    \draw [dashed] (6 ,5) -- (6,-5);
    \draw (-3,-5) -- (-3,5);
    \draw (-2,-5) -- (-2,5);
    \draw [black] (-4,-5) arc (0:180:3);
    \draw [black] (-5,-5) arc (0:180:2);
    \draw [black] (-6,-5) arc (0:180:1);

    \draw [black] (4,-5) arc (0:180:2);
    
    \draw [black] (5,5) arc (0:-180:3);
    \draw [black] (4,5) arc (0:-180:2);

    \draw [black] (-4,5) arc (0:-180:3);
    \draw [black] (-5,5) arc (0:-180:2);
    \end{tikzpicture}}} \; .
\end{equation*}

We will need the following result.
\begin{proposition}\label{prop: merging and moving caps}\cite[Proposition 3.9]{FK97}
Let $D$ and $D'$ be two TL diagrams that differ at the \textbf{left end} of the diagram by
\begin{equation*}
    D: \vcenter{\hbox{\begin{tikzpicture}[baseline={(current bounding box.center)}, scale=0.1]
    \draw [dashed] (-5,-5) -- (3,-5);
    \draw [dashed] (-5,5) -- (3,5);
    \draw [dashed] (-5,-5) -- (-5,5);
    \draw [black] (1,-5) arc (0:180:2);
    \draw [black] (-4,-5) -- (-4,5);
    \node at (0,-2) {y};
    \node at (-3,0) {x};
\end{tikzpicture}}} \quad
    D': \vcenter{\hbox{\begin{tikzpicture}[baseline={(current bounding box.center)}, scale=0.1]
    \draw [dashed] (-5,-5) -- (3,-5);
    \draw [dashed] (-5,5) -- (3,5);
    \draw [dashed] (-5,-5) -- (-5,5);
    \draw [black] (-0.5,-5) arc (0:180:2);
    \draw [black] (0,-5) -- (0,5);
    \node at (-2.5,-2) {y};
    \node at (1,0) {x};
\end{tikzpicture}}}
\end{equation*}
and are otherwise identical. Then $\coeffJW D = { x+y \brack x} \cdot \coeffJW D'$.
\end{proposition}
Given a TL diagram $D$, we will be interested in the following three operations.
\begin{itemize}
    \item[(R)] Rotate $D$ by one strand;
    \item[(M1)] Merge two neighboring nested caps (resp. cups) anywhere in $D$: 
    \begin{equation*}
    \vcenter{\hbox{\begin{tikzpicture}[baseline={(current bounding box.center)}, scale=0.1]
    \draw [dashed] (-5,-5) -- (5,-5);
    \draw [black] (-0.5,-5) arc (0:180:2);
    \draw [black] (4.5,-5) arc (0:180:2);
    \node at (-2.5,-2) {$x$};
    \node at (2.5,-2) {$y$};
\end{tikzpicture}}} \rightarrow
   \vcenter{\hbox{\begin{tikzpicture}[baseline={(current bounding box.center)}, scale=0.1]
    \draw [dashed] (-3,-5) -- (3,-5);
    \draw [black] (2,-5) arc (0:180:2);
    \node at (0,-2) {$x+y$};
\end{tikzpicture}}};
    \end{equation*}
    \item[(M2)] Move a nested cap through several strands at the left end of $D$:
    \begin{equation*}
    \vcenter{\hbox{\begin{tikzpicture}[baseline={(current bounding box.center)}, scale=0.1]
    \draw [dashed] (-5,-5) -- (3,-5);
    \draw [dashed] (-5,5) -- (3,5);
    \draw [dashed] (-5,-5) -- (-5,5);
    \draw [black] (1,-5) arc (0:180:2);
    \draw [black] (-4,-5) -- (-4,5);
    \node at (0,-2) {$y$};
    \node at (-3,0) {$x$};
\end{tikzpicture}}} 
    \rightarrow \vcenter{\hbox{\begin{tikzpicture}[baseline={(current bounding box.center)}, scale=0.1]
    \draw [dashed] (-5,-5) -- (3,-5);
    \draw [dashed] (-5,5) -- (3,5);
    \draw [dashed] (-5,-5) -- (-5,5);
    \draw [black] (-0.5,-5) arc (0:180:2);
    \draw [black] (0,-5) -- (0,5);
    \node at (-2.5,-2) {$y$};
    \node at (1,0) {$x$};
\end{tikzpicture}}}.
\end{equation*}
\end{itemize}

\begin{lemma}\label{lem: operations preserve degree and non-zero coefficient in Jones-Wenzl}
Let $D$ be a TL diagram with $\coeffJW D \neq 0$. Let $D'$ be another TL diagram that can be obtained from $D$ by a finite number of operations of the form (R), (M1) and (M2). Then $\coeffJW D' \neq 0$ and $\deg(D)  \in \{ \pm \deg(D') \}$.
\end{lemma}

\begin{proof}
If $D$ and $D'$ differ by an operation of the form (R), then $\coeffJW D' \neq 0$ by the rotatability of the Jones-Wenzl projector, and $\deg(D') = - \deg(D)$ by Lemma \ref{lem: rotation introduces minus on degree}. If $D$ and $D'$ differ by an operation of the form (M2), one can check using \eqref{eq: degree formula} that $\deg(D) = \deg(D')$ unless $x$ and $y$ are both odd. If $x$ is even, the parity of the caps don't change when moving them through the strands. If $y$ is even the total degree of a $y$-nested cap is $y/2 \cdot(g_s-f_t+g_t-f_s)$ no matter where it is in the diagram. When $x$ and $y$ are both odd, we have ${x+y \brack x} = 0$ by \eqref{quantumbinomclaim}. By Proposition \ref{prop: merging and moving caps} this implies $\coeffJW D  = {x+y \brack x} \coeffJW D' = 0$ which contradicts the assumption $\coeffJW D \neq 0$. Hence, $x$ and $y$ cannot be both odd if $\coeffJW D \neq 0$. Thus, we always have $\deg(D) = \deg(D')$ in the (M2) case and $0 \neq \coeffJW D  = {x+y \brack x} \coeffJW D'$ implies $ \coeffJW D' \neq 0$. Finally, (M1) can be obtained from (M2) and (R). More precisely, one can rotate (M2) by $x$ strands counterclockwise to obtain (M1) on the far left of the diagram. If the (M1) configuration is not on the far left, one can rotate until it is on the far left, merge the caps, and rotate back.
\end{proof}

We now verify the homogeneity of $\JW$.

\begin{proposition}\label{prop: Homogenity of Jones-Wenzl}
The Jones-Wenzl projector $\JW \in \TL_{m-1}(0)$ is homogeneous with $\deg(\JW) = 0$.
\end{proposition}

\begin{proof}
Consider a diagram of the form
\begin{equation} \label{hooray}
    D = \vcenter{\hbox{\begin{tikzpicture}[baseline={(current bounding box.center)}, scale=0.1]
    \draw [dashed] (-11,5) -- (1,5);
    \draw [dashed] (-11,-5) -- (1,-5);
    \draw [dashed] (-11,5) -- (-11,-5);
    \draw [dashed] (1,5) -- (1,-5);
    \draw (-5,-5) -- (-5,5);
        \node [font = \tiny]at (-8,-2) {$b$};
    \draw [black] (-6,-5) arc (0:180:2);
    \node [font = \tiny]at (-4,0) {$a$};
    \draw [black] (0,5) arc (0:-180:2);
    \node [font = \tiny]at (-1,2) {$b$};
\end{tikzpicture}}}
\end{equation}
for some $a,b$ with $a+ 2b= m-1$. Since $m-1$ is odd, $a$ is also odd. This implies that in $D$ we have $\# \{\text{even caps}\} =  \# \{ \text{odd cups} \}$ and $\# \{\text{odd caps} \} =  \# \{ \text{even cups}\}$. By \eqref{eq: degree formula}, this implies $\deg(D) = 0$.	
	
Now let $D$ be a TL diagram with $\coeffJW D \neq 0$. We need to show that $\deg(D) = 0$. By Lemma \ref{lem: operations preserve degree and non-zero coefficient in Jones-Wenzl} we can apply any number of operations of the form (R), (M1) or (M2) without changing the degree of $D$ (up to sign) or the condition $\coeffJW D \neq 0$. Our strategy is to apply these operations until we arrive at a diagram of the form \eqref{hooray}. This proves that $\JW$ is homogeneous of degree zero.

Suppose $D$ has $l-1$ through strands. Let $m = 2d$. Since $m-1$ is odd, $l-1$ is positive (and odd). Using (M1), we may merge neighboring nested cups/caps in $D$ until it is of the form
\begin{equation} \label{stdform}
    D = \vcenter{\hbox{\begin{tikzpicture}[baseline={(current bounding box.center)}, scale=0.1]
    \draw [dashed] (-11,5) -- (13,5);
    \draw [dashed] (-11,-5) -- (13,-5);
    \draw [dashed] (-11,5) -- (-11,-5);
    \draw [dashed] (13,5) -- (13,-5);
    \draw (-5,-5) -- (-5,5);
    \draw [black] (-6,-5) arc (0:180:2);
    \node [font = \scriptsize]at (-7,-2) {$b_1$};
    \draw [black] (-6,5) arc (0:-180:2);
    \node [font = \scriptsize]at (-7,2) {$a_1$};
    \draw (1,-5) -- (1,5);
    \draw [black] (0,-5) arc (0:180:2);
    \node [font = \scriptsize]at (-1,-2) {$b_2$};
    \draw [black] (0,5) arc (0:-180:2);
    \node [font = \scriptsize]at (-1,2) {$a_2$};
    \foreach \x in {3,4,5} {
    \node[circle, fill, inner sep=0.01cm] at (\x,0) {};
    }
    \draw (7,-5) -- (7,5);
    \draw [black] (12,-5) arc (0:180:2);
    \node [font = \scriptsize]at (11,-2) {$b_l$};
    \draw [black] (12,5) arc (0:-180:2);
    \node [font = \scriptsize]at (11,2) {$a_l$};
\end{tikzpicture}}}
\end{equation}
for some $a_1, \ldots, a_l, b_1, \ldots, b_l \ge 0$. By counting boundary strands we have $2(a_1 + ... + a_l) = 2(b_1 + ... + b_l) = m-l$. In particular, $2b_1$ is at most $m - 2$, so $b_1 \le d-1$. We work case by case.

{\bf Case 1}: $b_1 < d-1$ and $a_1  + 2b_1 \le  m-1$. We then distinguish two subcases.

{\bf Case 1.1}: $a_1> 0$.

Since $a_1 + 2b_1 \le m-1$, rotation by $a_1$ many strands (counterclockwise) results in
\begin{equation*}
     \vcenter{\hbox{\begin{tikzpicture}[baseline={(current bounding box.center)}, scale=0.1]
    \draw [dashed] (-5,-5) -- (7,-5);
    \draw [dashed] (-5,5) -- (7,5);
    \draw [dashed] (-5,-5) -- (-5,5);
    \draw [black] (1,-5) arc (0:180:2);
    \draw [black] (-4,-5) -- (-4,5);
    \node[font= \tiny] at (0,-2) {$b_1$};
    \node[font= \tiny] at (-2,1) {$a_1$};
        \foreach \x in {4,5,6} {
    \node[circle, fill, inner sep=0.01cm] at (\x,0) {};
    }
\end{tikzpicture}}},
\end{equation*}
that is, none of the strands in the $b_1$ nested caps are moved to the top by this rotation. Then use (M2) to obtain
\begin{equation*}
     \vcenter{\hbox{\begin{tikzpicture}[baseline={(current bounding box.center)}, scale=0.1]
    \draw [dashed] (-5,-5) -- (7,-5);
    \draw [dashed] (-5,5) -- (7,5);
    \draw [dashed] (-5,-5) -- (-5,5);
    \draw [black] (0,-5) arc (0:180:2);
    \draw [black] (1,-5) -- (1,5);
    \node[font= \tiny] at (-1,-2) {$b_1$};
    \node[font= \tiny] at (-0.5,1) {$a_1$};
        \foreach \x in {4,5,6} {
    \node[circle, fill, inner sep=0.01cm] at (\x,0) {};
    }
\end{tikzpicture}}}.
\end{equation*}
Now, after possibly merging neighboring nested cups and and caps to obtain a diagram of the form \eqref{stdform}, we are in Case 1.2.

{\bf Case 1.2}: $a_1 = 0$.

In this case, $D$ is of the form
\begin{equation*}
    D = \vcenter{\hbox{\begin{tikzpicture}[baseline={(current bounding box.center)}, scale=0.1]
    \draw [dashed] (-11,5) -- (-1,5);
    \draw [dashed] (-11,-5) -- (-1,-5);
    \draw [dashed] (-11,5) --(-11,-5);
    \draw (-5,-5) -- (-5,5);
    \draw [black] (-6,-5) arc (0:180:2);
    \node [font = \scriptsize]at (-7,-2) {$b_1$};

    \foreach \x in {-4,-3,-2} {
    \node[circle, fill, inner sep=0.01cm] at (\x,0) {};
    }
\end{tikzpicture}}}.
\end{equation*}
Rotation by one strand counterclockwise results in
\begin{equation*}
     \vcenter{\hbox{\begin{tikzpicture}[baseline={(current bounding box.center)}, scale=0.1]
    \draw [dashed] (-14,5) -- (0,5);
    \draw [dashed] (-14,-5) -- (0,-5);
    \draw [dashed] (-14,5) -- (-14,-5);
    \draw [black] (-5,-5) arc (0:180:4);
    \draw [black] (-7,-5) arc (0:180:2);
    \node [font = \scriptsize]at (-8,-2) {$b_1$};
    \foreach \x in {-3,-2,-1} {
    \node[circle, fill, inner sep=0.01cm] at (\x,0) {};
    }
\end{tikzpicture}}},
\end{equation*}
that is, the condition $b_1 < d-1$ implies that the bottom endpoint of the first thru-strand was not moved to the top by this rotation.
Merging all neighboring nested cups/caps, we obtain
\begin{equation*}
D' =\vcenter{\hbox{\begin{tikzpicture}[baseline={(current bounding box.center)}, scale=0.1]
    \draw [dashed] (-11,5) -- (14,5);
    \draw [dashed] (-11,-5) -- (14,-5);
    \draw [dashed] (-11,5) -- (-11,-5);
    \draw [dashed] (14,5) -- (14,-5);
    \draw (-5,-5) -- (-5,5);
    \draw [black] (-6,-5) arc (0:180:2);
    \node [font = \scriptsize]at (-7,-2) {$b'_1$};
    \draw [black] (-6,5) arc (0:-180:2);
    \node [font = \scriptsize]at (-7,2) {$a'_1$};
    \draw (1,-5) -- (1,5);
    \draw [black] (0,-5) arc (0:180:2);
    \node [font = \scriptsize]at (-1,-2) {$b'_2$};
    \draw [black] (0,5) arc (0:-180:2);
    \node [font = \scriptsize]at (-1,2) {$a'_2$};
    \foreach \x in {3,4,5} {
    \node[circle, fill, inner sep=0.01cm] at (\x,0) {};
    }
    \draw (7,-5) -- (7,5);
    \draw [black] (12,-5) arc (0:180:2);
    \node [font = \scriptsize]at (12,-2) {$b'_{l'}$};
    \draw [black] (12,5) arc (0:-180:2);
    \node [font = \scriptsize]at (12,2) {$a'_{l'}$};
\end{tikzpicture}}}
\end{equation*}
where $b'_1 \ge b_1 + 1$.

Thus if we are in Case 1, we can replace the original diagram $D$ with another diagram $D'$ also in the form \eqref{stdform}, where $b'_1$ is larger than $b_1$. We can repeat this process until we are no longer in Case 1.


{\bf Case 2}: $b_1 = d-1$ but $a_1  + 2b_1 \le  m-1$.

Since $2b_1+1 = m-1$ is the total number of boundary points on bottom, there is only one through strand. Note that $0 \le a_1 \le 1$. If $a_1 = 0$ then
\begin{equation*}
    D = \vcenter{\hbox{\begin{tikzpicture}[baseline={(current bounding box.center)}, scale=0.1]
    \draw [dashed] (-11,5) -- (1,5);
    \draw [dashed] (-11,-5) -- (1,-5);
    \draw [dashed] (-11,5) --(-11,-5);
    \draw [dashed] (1,5) -- (1,-5);
    \draw (-5,-5) -- (-5,5);
    \draw [black] (-6,-5) arc (0:180:2);
    \node [font = \tiny]at (-8,-2) {d-1};
    \draw [black] (0,5) arc (0:-180:2);
    \node [font = \tiny]at (-1,2) {d-1};
\end{tikzpicture}}}.
\end{equation*}
If $a_1 = 1$ then $D$ is of the form 
\begin{equation*}
    D = \vcenter{\hbox{\begin{tikzpicture}[baseline={(current bounding box.center)}, scale=0.1]
    \draw [dashed] (-11,5) -- (1,5);
    \draw [dashed] (-11,-5) -- (1,-5);
    \draw [dashed] (-11,5) --(-11,-5);
    \draw [dashed] (1,5) -- (1,-5);
    \draw (-5,-5) -- (-5,5);
    \draw [black] (-6,5) arc (0:-180:2);
    \draw [black] (-6,-5) arc (0:180:2);
    \node [font = \tiny]at (-8,-2) {d-1};
    \draw [black] (0,5) arc (0:-180:2);
    \node [font = \tiny]at (-1,2) {d-2};
\end{tikzpicture}}}.
\end{equation*}
We can then rotate counterclockwise by one strand to obtain
\begin{equation*}
    \vcenter{\hbox{\begin{tikzpicture}[baseline={(current bounding box.center)}, scale=0.1]
    \draw [dashed] (-5,-5) -- (5,-5);
    \draw [dashed] (-5,5) -- (5,5);
    \draw [dashed] (5,5) -- (5,-5);
    \draw [dashed] (-5,-5) -- (-5,5);
    \draw [black] (3,-5) arc (0:180:2);
    \draw [black] (3,5) arc (0:-180:2);
    \draw [black] (-4,-5) -- (-4,5);
    \node[font= \tiny] at (1,-2) {d-1};
    \node[font= \tiny] at (1,2) {d-1};
\end{tikzpicture}}}
\end{equation*}
and apply (M2) to obtain again 
\begin{equation*}
    \vcenter{\hbox{\begin{tikzpicture}[baseline={(current bounding box.center)}, scale=0.1]
    \draw [dashed] (-11,5) -- (1,5);
    \draw [dashed] (-11,-5) -- (1,-5);
    \draw [dashed] (-11,5) --(-11,-5);
    \draw [dashed] (1,5) -- (1,-5);
    \draw (-5,-5) -- (-5,5);
    \draw [black] (-6,-5) arc (0:180:2);
    \node [font = \tiny]at (-8,-2) {d-1};
    \draw [black] (0,5) arc (0:-180:2);
    \node [font = \tiny]at (-1,2) {d-1};
\end{tikzpicture}}}.
\end{equation*}
Either way this is in the desired form from \eqref{hooray}.

{\bf Case 3}: $a_1  + 2b_1 > m-1$.

In this case, we apply a few more operations. Let $k = m-1-2b_1$, noting that $k < a_1$. First rotate by $k$ strands (counterclockwise) to obtain
    \begin{equation*}
    \vcenter{\hbox{\begin{tikzpicture}[baseline={(current bounding box.center)}, scale=0.1]
    \draw [dashed] (-12,5) -- (7,5);
    \draw [dashed] (-12,-5) -- (7,-5);
    \draw [dashed] (-12,5) -- (-12,-5);
    \draw [dashed] (7,5) -- (7,-5);
    \draw (-5,-5) -- (-5,5);
        \node [font = \scriptsize]at (-4,0) {$k$};
    \draw [black] (-6,5) arc (0:-180:2);
    \node [font = \scriptsize]at (-9,2) {$a_1$-$k$};
    \draw [black] (1,-5) arc (0:180:2);
    \node [font = \scriptsize]at (0,-2) {$b_1$};
      \draw [black] (6,5) arc (0:-180:4);
          \draw [black] (0,2) rectangle (4,5);
    \node[font=\scriptsize] at (2,3.2) {\textbf{*}};
\end{tikzpicture}}}.
\end{equation*}
The starred diagram in the top right corner consists of cups which can be merged with (M1) or nested, producing the diagram
    \begin{equation*}
     D' = \vcenter{\hbox{\begin{tikzpicture}[baseline={(current bounding box.center)}, scale=0.1]
    \draw [dashed] (-12,5) -- (6,5);
    \draw [dashed] (-12,-5) -- (6,-5);
    \draw [dashed] (-12,5) -- (-12,-5);
    \draw [dashed] (6,5) -- (6,-5);
    \draw (-5,-5) -- (-5,5);
    \node [font = \scriptsize]at (-4,0) {$k$};
    \draw [black] (-6,5) arc (0:-180:2);
    \node [font = \scriptsize]at (-9,2) {$a_1$-$k$};
    \draw [black] (2,-5) arc (0:180:2);
    \node [font = \scriptsize]at (1,-2) {$b_1$};
    \draw [black] (2,5) arc (0:-180:2);
    \node [font = \scriptsize]at (2,2) {$x_3$};
\end{tikzpicture}}}
\end{equation*}
where $x_3 = m - 1 - a_1 - b_1$.
Let us write $x_1 = a_1$ and $x_2 = b_1$. Up to rotation, we can thus depict $D'$ as a diagram
\begin{equation*}
D' = \vcenter{\hbox{\begin{tikzpicture}[baseline={(current bounding box.center)}, scale=0.7]
  \draw[dashed] (0,0) circle (1);
  \draw [black] (0.3,0.95) arc (0:-180:0.3);
  \node [font = \scriptsize]at (0,0.4) {$x_1$};
  \draw [black] (-1,-0.2) arc (120:-60:0.3);
    \node [font = \scriptsize]at (-0.3,-0.2) {$x_2$};
  \draw [black] (1,-0.2) arc (60:240:0.3);
  \node [font = \scriptsize]at (0.3,-0.2) {$x_3$};
\end{tikzpicture}}}
\end{equation*}
with
\begin{equation}\label{eq: the xi sum to m-1}
    x_1+x_2+x_3 = m-1.
\end{equation}
If $x_1 + 2x_2, x_2 + 2x_3$ and $x_3 + 2x_1$ are all greater than $m-1$ we get
\begin{align*}
    3(m-1) & = 3(x_1 + x_2 + x_3) \\
    &= (x_1 + 2x_2) + (x_2 + 2x_3) + (x_3 + 2x_1)\\
    & >  3(m-1)
\end{align*}
which is a contradiction. Hence, without loss of generality, we may assume
\begin{equation}\label{eq: wlog x2 +2x1 is at most m-1}
    x_2 + 2x_3 \le m-1
\end{equation}
(the other cases follow from this by rotation and reindexing). Then, after an appropriate rotation, we get
    \begin{equation*}
     \vcenter{\hbox{\begin{tikzpicture}[baseline={(current bounding box.center)}, scale=0.1]
    \draw [dashed] (-7,5) -- (8,5);
    \draw [dashed] (-7,-5) -- (8,-5);
    \draw [dashed] (-7,5) -- (-7,-5);
    \draw [dashed] (8,5) -- (8,-5);
    \draw (-5,-5) -- (-5,5);
    \node [font = \tiny]at (-3,1) {$x_2$};
    \draw [black] (0,-5) arc (0:180:2);
    \node [font = \tiny]at (-1,-2) {$x_3$};
    \draw (1,-5) -- (1,5);
    \node [font = \tiny]at (3,0) {$k'$};
    \draw [black] (6,5) arc (0:-180:2);
    \node [font = \tiny]at (6,2) {$x_3$};
\end{tikzpicture}}}
\end{equation*}
where $k' = x_1-x_3$ (note that \eqref{eq: the xi sum to m-1} and \eqref{eq: wlog x2 +2x1 is at most m-1} imply that $x_1 \ge x_3$). Finally, using (M2) one more time, we get
    \begin{equation*}
     \vcenter{\hbox{\begin{tikzpicture}[baseline={(current bounding box.center)}, scale=0.1]
    \draw [dashed] (-11,5) -- (3,5);
    \draw [dashed] (-11,-5) -- (3,-5);
    \draw [dashed] (-11,5) -- (-11,-5);
    \draw [dashed] (3,5) -- (3,-5);
    \draw (-5,-5) -- (-5,5);
        \node [font = \tiny]at (-1,-1) {$x_2$+$k'$};
    \draw [black] (-6,-5) arc (0:180:2);
    \node [font = \tiny]at (-8,-2) {$x_3$};
    \draw [black] (0,5) arc (0:-180:2);
    \node [font = \tiny]at (-1,2) {$x_3$};
\end{tikzpicture}}},
\end{equation*}
which has the desired form \eqref{hooray}.
\end{proof}

Since $\JW$ is homogeneous of degree $0$, we get that the corresponding Jones-Wenzl relation \cite[(10.7i)]{EMTW} is also homogeneous. Together with Lemma \ref{lemma: homogeneity of JW in the easy case} this verifies the homogeneity of the Jones-Wenzl relation in all cases. Thus, by Proposition \ref{prop: gradings extend if and only JW is homogeneous}, we obtain the following universal grading on the Hecke category.

\begin{theorem}\label{thm:mostgeneral}
Let $V$ be a realization with a $W$-invariant grading valued in an abelian group $\Gamma$ for which the simple roots are homogeneous. Let $\Lambda_{\sim}$ be the free abelian group on generators $\{ f_s, g_s\}_{s \in S}$ modulo the relation $f_s + g_s = f_t + g_t$ whenever $\langle \alpha_s, \alpha_t^{\vee} \rangle \neq 0$. Then the grading on $V \subset R$ extends to a grading on $\HC$ by the group $(\Gamma \times \Lambda_{\sim})/I$ where $I$ is the submodule generated by $(-\deg(\alpha_s),f_s + g_s)$. The generating morphisms have degrees as in \eqref{eq:multigrading}. Moreover, this is the universal monoidal grading which refines the original grading and extends the grading on $V$.
\end{theorem}

\subsection{Grading based on characteristic}

As an important special case, we obtain the $p$-graded Hecke category. We begin with a choice of prime number $p$, the characteristic of the base field; then we choose a realization so
that $\langle \alpha_s^\vee, \alpha_t \rangle = 0$ whenever $m_{st} = 2p^k$ for some $k$.

\begin{defn} Suppose that $\Bbbk$ has characteristic $p$, possibly zero. Recall the equivalence relation $\sim_p$ from Definition \ref{defn:simp}. For each equivalence class $J \subset S$ choose a realization $V_J$ of the parabolic subgroup $W_J = \langle s \in J \rangle$. Extend $V_J$ to a representation of $W$ for which all simple reflections $t \notin J$ act trivially. Then let
        \begin{equation*}
            V = \bigoplus_{J \in S / \sim_p} V_J.
        \end{equation*}
In the following lemma we verify that $V$ is a realization of $W$, with the natural choice of simple roots and coroots extended from each $V_J$. Any realization constructed in such a way is called \emph{$p$-adapted}. \end{defn}

\begin{lemma} \label{lem:coloradapted} In the previous definition, $V$ is a well-defined realization of $W$, for which the Hecke category is well-defined. \end{lemma}
	
\begin{proof} First we confirm that $V$ is a representation of $W$. The braid relations within each given parabolic subgroup $W_J$ are immediate, as are the quadratic relations. If $s
\not\sim_p t$ then $s$ and $t$ commute on $V$ by construction. Moreover, if $m_{st}$ is finite then $m_{st} = 2p^k$ which is even, so the braid relation holds.

All we need to check for $V = \bigoplus_{J \in S/\sim_p} V_i$ to be a realization in the sense of \cite{HaziRotatable} is that for $s,t$ with $m_{st} < \infty$ the corresponding Jones-Wenzl projector exists and is rotatable. If $s \sim_p t$ this follows since each $V_J$ is a realization. If $s \not \sim_p t$ we have $m_{st} =2p^k$ and the result follows by Lemma \ref{lemma: m is 2 p^k}. \end{proof}

For any $J \in S / \sim_p$, we let $e_J$ denote the element of $\Z^{S / \sim_p}$ which is $1$ in the $J$-th factor and zero elsewhere. Then $\{e_J\}$ is a $\Z$-basis for this abelian group. For any $s \in S$, write $[s] \in S/\sim_p$ for its equivalence class.

\begin{corollary}\label{corollary: color graded Hecke category} Let $V$ be a $p$-adapted realization of $(W,S)$ over $\Bbbk$. Equip $V$ with a grading by $\Z^{S / \sim_p}$ such that $V_J$ has degree $2e_J$ for each $J \in S/\sim_p$. This extends to a grading of the Hecke category by setting
     \begin{subequations} 
     \begin{equation*}  \deg \left( \startdotred \right) = \deg \left( \finaldotred \right) = e_{[s]}, \qquad \deg \left( \splitred \right) = 
   \deg \left( \mergered \right) = -e_{[s]}, \end{equation*}
     \begin{equation*}  \deg \left( \sbotttop \right) =0. \end{equation*}
     \end{subequations}
Indeed, given a $p$-adapted realization, this grading of $\HC$ is universal under the constraints that each $V_J$ is in a single degree, and that the degree is invariant under vertical symmetry.
\end{corollary}

\begin{proof} This follows directly from Theorem \ref{thm:mostgeneral}. \end{proof}

\section{The Grothendieck group} \label{sec:catfn}

The goal of this section is to compute the Grothendieck group of the Hecke category when equipped with a more general grading.

We do \emph{not} assume that $A$ refines the standard $\Z$-grading! This flexibility will be necessary to be able to categorify unequal parameter Hecke algebras where all
parameters are powers of a single variable $v$. For example, one might have $A = \Z$ with $f_s = g_s = 1$ and $f_t = g_t = 2$ for some $s \ne t$; it is impossible to specialize this
grading back to the original grading.

We still assume that the generating diagrams of $\HC$ are homogeneous.

\subsection{Size of morphism spaces}

\begin{assumption} \label{assume} Let $A$ be an abelian group. Let $V$ have a good $A$-grading, and assume it extends to an $A$-grading on $\HC$ where the generating diagrams are homogeneous. We define $f_s, g_s \in A$ for each $s \in S$ so that \eqref{eq:onecolor} holds. We assume that $A$ is equipped with a \emph{bar involution}, which
satisfies $\bar{f}_s = -g_s$ and $\bar{g}_s = -f_s$, and $\bar{a} = -a$ whenever $a$ is the degree of a homogeneous element element of $V$. We also assume that there exists a group
homorphism $\grade \colon A \to \Z$ such that all degrees of homogeneous elements of $V$ are sent to $\Z_{> 0}$, though we do not fix such a homomorphism. \end{assumption}

\begin{assumption} \label{assume2} We further assume that $f_s = f_t$ and $g_s = g_t$ whenever $m_{st}$ is odd. \end{assumption}
	
The last condition of Assumption \ref{assume} implies that, under the $\Z$-grading induced by $\grade$, the ring $R=\End(\emptyset) = S(V)$ is non-negatively graded with $R^0 = \Bbbk$. It also implies
that $\grade(f_s+g_s) > 0$ for all $s$, since $\deg(\alpha_s) = f_s + g_s$. We use this positivity crucially in Lemma \ref{lemma: basis of split grothendieck group of H}, though we wonder whether it is
possible to bypass this argument.

Removing Assumption \ref{assume2} introduces technicalities which we postpone until \S\ref{subsec:rescaled}. 

\begin{lemma} Under Assumption \ref{assume}, the $2m_{s,t}$-valent vertex obeys \eqref{eq:2mvalent}. Under Assumption \ref{assume2}, it has degree zero. \end{lemma}
	
\begin{proof} The proof of Proposition \ref{prop: gradings extend if and only JW is homogeneous} applies using our assumptions, as noted therein. \end{proof}

The bar involution is introduced so that the following lemma holds.

\begin{lemma} \label{lem:barandflip} Under Assumption \ref{assume}, flipping a diagram upside down induces an anti-automorphism of $\HC$ which acts on the degree by minus the bar involution. \end{lemma}

\begin{proof} We need only prove this statement for the generating morphisms, where it is straightforward from \eqref{eq:onecolor} and the previous lemma. \end{proof}

We let $B_s$ denote the generating object of $\HC$, and for a sequence $\un{w}$ of simple reflections we write $B_{\un{w}}$ for the corresponding Bott-Samelson object (the tensor product of various $B_s$). We denote by $\HC_A$ the category with objects $X(a)$, where $X$ is a Bott-Samelson object and $a \in A$, and with morphism spaces
\begin{equation*}
    \Hom_{\HC_A}(X(a), Y(b)) = \Hom_{\HC}^{b-a}(X,Y).
\end{equation*}
We denote by $\Kar (\HC_A)$ the Karoubi envelope of the additive closure of $\HC_A$ and we denote the corresponding split Grothendieck group by
\begin{equation*}
    [\HC]_A = [\Kar (\HC_A)]_{\oplus}.
\end{equation*}
Note that $[\HC]_A$ has a canonical $\mathbb{Z}[A]$-module structure where $v^a$ acts by the shift $(-a)$.

The choice of grading does not affect the fact that morphism spaces in $\HC$ are finite-rank free modules over $R$, which is an $A$-graded ring. Pick some $\grade$ as in Assumption \ref{assume}. In the induced $\Z$-grading, $R$ is positively graded and finite dimensional in each degree. Finite generation over $R$ implies that morphism spaces in the induced $\Z$-grading on $\HC$ are finite-dimensional in each degree, and hence the same holds in the original $A$-grading. Consequently, the Grothendieck group $[\HC]_A$ has a $\Z[A]$-sesquilinear Hom pairing defined by
\[ ([M],[N]) := grk_R \Hom_{\HC}^{\bullet} (M,N)  \]
where for any a graded free $R$-module $M = R(a_1) \oplus ... \oplus R(a_k)$ we write \[ grk_R M := v^{-a_1} + ... + v^{-a_k} .\]

What makes computing the Grothendieck group subtle is that the classification of indecomposable objects in the original grading on $\HC$ does not immediately yield a classification of indecomposable objects in $A$-graded $\HC$,
see the discussion in \S\ref{intro-grothgroup}. Instead, we will take the same technology which was used in \cite{EWGr4sb} to classify the indecomposable objects, and argue that it still applies to
$A$-graded $\HC$: namely, the technology of light leaves, double leaves, and object-adapted cellular categories.

%
%
%

Let us briefly recall the construction of light leaves and double leaves from \cite[Construction 6.1]{EWGr4sb}. Pick for each $w \in W$ a reduced expression $\underline{w}$. Then for any $w,x \in W$ and any subexpression $\underline{e} \subset \underline{w}$ expressing $x$, one can define a distinguished Soergel diagram $LL_{\underline{w}, \underline{e}} \in \Hom (B_{\underline{w}}, B_{\underline{x}})$ which we call a \emph{light leaf}\footnote{The light leaf depends on further choices, but we make these choices once and for all and omit them from the notation. We have chosen exactly one diagram $LL_{\underline{w}, \underline{e}}$ for each pair $(\un{w}, \un{e})$.}. Given subexpressions $\underline{e} \subset \underline{w}$ and $\underline{f} \subset \underline{y}$ such that $\underline{e}$ and $\underline{f}$ both express $x$, one defines the corresponding \emph{double leaf}
\begin{equation*}
    \mathbb{LL}^x_{\underline{f}, \underline{e}}= \overline{LL}_{\underline{y}, \underline{f}} \circ LL_{\underline{w} ,\underline{e}} \in \Hom(B_{\underline{w}}, B_{\underline{y}})
\end{equation*}
where $\overline{LL}_{\underline{y}, \underline{f}}$ is given by flipping $LL_{\underline{y}, \underline{f}}$ upside down. The set of all $\mathbb{LL}^x_{\underline{f}, \underline{e}}$ forms a basis of $\Hom(B_{\underline{w}}, B_{\underline{y}})$ as a free (left) $R$-module, known as the \emph{double leaves basis}. Note that double leaves are diagrams, so they are homogeneous in the $A$-grading, and form a homogeneous basis for morphism spaces as $A$-graded $R$-modules. This basis equips $\HC_A$ with the structure of an object adapted cellular category (see \cite[Definition 11.29]{EMTW} or \cite[Definition 2.4]{ELauda}).

\begin{lemma}\label{lemma: basis of split grothendieck group of H}
Pick for each $w \in W$ a reduced expression $\underline{w}$, as we did in the construction of double leaves. Under Assumption \ref{assume}, the set $\{[B_{\un{w}}]\}_{w \in W}$ forms a $\mathbb{Z}[A]$-basis of $[\HC]_A$. Moreover, each $B_{\underline{w}} \in  \Kar( \HC_A)$ has a unique indecomposable summand (abusively denoted by $B_w$) which does not appear as a direct summand of $B_{\underline{y}}$ for any $y < w$, and this gives a complete set of indecomposable objects up to grading shift and isomorphism. \end{lemma}
	
        

This proof is analogous to the one in \cite[Theorem 6.25]{EWGr4sb}, which was generalized by \cite[Proposition 2.24]{ELauda}. The proof was rewritten (ahistorically) in \cite[Theorem 11.39]{EMTW} to use the technology of \cite{ELauda}, and this is the version which we follow. 

\begin{proof} Note that $\HC$ is an object adapted cellular category in the sense of \cite[Definition 11.29]{EMTW}. The cellular basis is the double leaves basis which is homogeneous with respect to
the $A$-grading since it consists of Soergel diagrams. The proof of \cite[Theorem 11.39]{EMTW} then applies word by word after replacing all $\mathbb{Z}$-degrees with $A$-degrees. At some point, the
proof uses that any degree $0$ element in a certain graded $R$-module $J_{\underline{w}}(X) \subset \End_{\underline{w}}(X) / \End_{\not \ge \underline{w}}(X)$ is nilpotent. For the original
$\mathbb{Z}$-grading this is proved in \cite[Lemma 11.41]{EMTW} under the assumption that $R =\End(\emptyset)$ is non-negatively graded and $R^0 = \Bbbk$. This result still holds for an $A$-grading
under Assumption \ref{assume}. Any element with $A$-degree $0$ would still have degree $0$ in the induced $\Z$-grading under $\grade$, thus must be nilpotent. \end{proof}

As a consequence of Lemma \ref{lemma: basis of split grothendieck group of H} we know the size of $[\HC]_A$, and as a consequence of the next lemma we know the Hom form on a basis
(the basis of symbols of certain Bott-Samelsons, not the basis of symbols of indecomposables).

\begin{defn} \label{Def:defect} Let $\eb \subset \un{w}$ be a subexpression, and label each index of $\eb$ with $U0$, $U1$, $D0$, or $D1$ as in \cite[\S 2.4]{EWGr4sb}. Then the \emph{defect} of $\eb$, denoted $\defect(\eb)$, is the element of $A$ obtained by summing up the defect of each index: if the $i$-th index in $\un{w}$ is $s$, then the contribution of the $i$-th term in $\eb$ is
\begin{equation} \defect(U1) = 0, \quad \defect(U0) = f_s, \quad \defect(D0) = -g_s, \quad \defect(D1) = f_s - g_s. \end{equation}
\end{defn}

\begin{lemma} \label{lem:degLL} Under Assumptions \ref{assume} and \ref{assume2}, in the $A$-grading on $\HC$ we have
\begin{equation} \deg(LL_{\un{w}, \un{e}}) = -\overline{\defect(\un{e})},
	 \qquad \deg(\overline{LL}_{\un{w}, \un{e}}) = \defect(\un{e}). \end{equation}
\end{lemma}

\begin{proof} We leave the reader to review \cite[Construction 6.1]{EWGr4sb}, and make only the following observations. Light leaves are built from: rex moves (degree zero), an enddot for each U0, an identity map for each U1, a merging trivalent vertex for each D0, and a cap for each D1. Thus it is easy to verify the first equality, and Lemma \ref{lem:barandflip} yields the second. \end{proof}
	
\begin{remark} Without Assumption \ref{assume2} the degrees of light leaves are more complicated to compute, since rex moves do not have degree zero. We address this in Remark \ref{rmk:degrees of light leaves}. \end{remark}

\subsection{Hecke algebras with unequal parameters and categorification} \label{subsec:heckeunequal}

Our goal will be to match $[\HC]_A$ with a Hecke algebra with unequal parameters $H_A$. The standard literature on Hecke algebras with unequal parameters is only adapted to the case where $f_s = g_s$, so we point out that assumption now. We remove the assumption in the next section.

\begin{assumption} \label{assume3} We further assume that $f_s = g_s$ for all $s \in S$. \end{assumption}

Note that Assumption \ref{assume3} implies Assumption \ref{assume2} in the absence of 2-torsion, since $f_s + g_s = f_t + g_t$ whenever $m_{st}$ is odd.

\begin{defn} Let $(W,S)$ be a Coxeter system. We place an equivalence relation $\sim_{\uneq}$ on $S$ where $s \sim_{\uneq} t$ if and only if $s$ and $t$ are conjugate in $W$. It is well known that
$\sim_{\uneq}$ is generated by \begin{equation} s \sim_{\uneq} t \text{ if } m_{st} \text{ is odd (and finite).} \end{equation} \end{defn}


\begin{defn} Let $(W,S)$ be a Coxeter system. Let $A_{\uneq} = \Z^{S/\sim_{\uneq}}$. Its group algebra is generated by invertible elements $v_s := v^{1_s}$ for each $s \in S$, modulo the relation $v_s = v_t$ when $s \sim_{\uneq} t$. The \emph{Hecke algebra with formal unequal parameters} $H_{\uneq}$ is the $\Z[A_{\uneq}]$-algebra with generators $\delta_s$ ($s \in S$) modulo the braid relations and the 
	quadratic relation
    \begin{equation} \label{unequalquadratic} 
    \delta_s^2 = (v_s^{-1} - v_s)\delta_s + 1.
    \end{equation}
Equivalently, we have
\begin{equation} (\delta_s - v_s^{-1})(\delta_s + v_s) = 0. \end{equation}
\end{defn}

\begin{defn} For any abelian group $A$, one can consider the group algebra $\Z[A]$. An \emph{($A$-valued) parameter function} is a group homomorphism $L : A_{\uneq} \rightarrow A$. A parameter function gives rise to a ring homomorphism $\Z[A_{\uneq}] \to \Z[A]$. Let $H_L$ denote the corresponding specialization of $H_{\uneq}$. If $L$ is understood, we also write $H_A$ for $H_L$.
\end{defn}

For example, the usual Hecke algebra is the specialization $H_\Z$ for the parameter function where $L(1_s) = 1$ for all $s$, where the group algebra of $\Z$ is identified with $\Z[v,v^{-1}]$. A traditional \emph{Hecke algebra with unequal parameters} is a specialization $H_L$ for some $\Z$-valued parameter function $L$. Everything said below about $H_{\uneq}$ applies equally to $H_L$ for any parameter function.

By \cite[Theorem 7.1]{HumpCox}, $H_{\uneq}$ is a free $\Z[A_{\uneq}]$-module with basis $\delta_w$ for $w \in W$, where $\delta_w$ is a product of various $\delta_s$ along a reduced expression for $w$. Let
\begin{equation} b_s := \delta_s + v_s. \end{equation}
It is well-known how to deduce that
\begin{equation}\label{dwbs}
	\delta_w b_s = \begin{cases} \delta_{ws} + v_s \delta_w, & \text{ if } ws > w,\\
\delta_{ws} + v_s^{-1} \delta_w, & \text{ if } ws < w. \end{cases} \end{equation}

Let $\un{w} (s_1, \ldots, s_d)$ be an expression in $S$. The product $b_{\un{w}} := b_{s_1} b_{s_2} \cdots b_{s_d}$ can be expanded in the basis $\{\delta_w\}$ by iterating the formula \eqref{dwbs}. For ordinary Hecke algebras the resulting expression is called the Deodhar defect formula, and it adapts easily.

\begin{defn} \label{Def:defect2} Let $\eb \subset \un{w}$ be a subexpression, and label each index of $\eb$ with $U0$, $U1$, $D0$, or $D1$ as in \cite[\S 2.4]{EWGr4sb}. Then the \emph{defect} of $\eb$, denoted $\defect(\eb)$, is the element of $A_{\uneq}$ obtained by summing up the defect of each index: if the $i$-th index in $\un{w}$ is $s$, then the contribution of the $i$-th term in $\eb$ is
\begin{equation} \defect(U1) = 0, \quad \defect(U0) = 1_s, \quad \defect(D0) = -1_s, \quad \defect(D1) = 0. \end{equation}
\end{defn}

\begin{lemma}\label{lem:defect} (c.f. Deodhar's defect formula \cite[Lemma 2.10]{EWGr4sb}) We have
	\begin{equation} b_{\un{w}} = \sum_{\eb \subset \un{w}} v^{\defect(\eb)} \delta_{\un{w}^{\eb}}, \end{equation}
where $\un{w}^{\eb}$ is the element expressed by the subexpression, see \cite[\S 2.4]{EWGr4sb}. \end{lemma}

\begin{proof} The proof is straightforward and identical to the classical case, mutatis mutandis. \end{proof}

Here is our categorification theorem.

\begin{theorem}\label{thm:categorification theorem}
    Under Assumptions \ref{assume} and \ref{assume2} and \ref{assume3}, the map $b_s \mapsto [B_s]$ induces an isomorphism of $\mathbb{Z}[A]$-algebras
    \begin{equation*}
        H_{A} \overset{\sim}{\rightarrow }[\HC]_A.
    \end{equation*}
	Here we use the $A$-valued parameter function where $1_s \mapsto f_s = g_s$. Under this isomorphism, the Hom form on $[\HC]_A$ agrees with the standard sesquilinear form on $H_A$.
\end{theorem}

\begin{proof}
Choose for each $w \in W$ a reduced expression $\underline{w}$ and consider the diagrammatic character $ch: [\HC]_A \rightarrow H_A$ defined by
\begin{equation}
    ch: [B] \mapsto \sum_{w \in W} grk_{R} \Hom_{\ge w}(B_{\underline{w}}, B)\delta_w.
\end{equation}
Here $\Hom_{\ge w}$ is the morphism space in the quotient category obtained by killing all double leaves that factor through some $y \not \ge w$. The map $ch$ is clearly a morphism of $\mathbb{Z}[A]$-modules. Note that $\Hom_{\ge w}(B_{\underline{w}},B_{\underline{x}})$ has the basis $\{\overline{LL}_{\underline{x}, \underline{f}} \}$ ranging over subexpressions $\un{f} \subset \un{x}$ which express $w$. It now follows from Lemma \ref{lem:defect}  and Lemma \ref{lem:degLL} that $ch([B_{\underline{x}}]) = b_{\underline{x}}$. In particular, $ch$ is multiplicative under concatenation of sequences (the monoidal product). It follows from Lemma \ref{lemma: basis of split grothendieck group of H} that the $[B_{\underline{x}}]$ span $[\HC]_A$. This shows that $ch$ is multiplicative and thus an algebra homomorphism. Moreover, choosing a reduced expression for each $w \in W$ we obtain a basis $\{[B_{\underline{w}}]\}$ for $[\HC]_A$ over $\Z[A]$. After applying $ch$ to this basis, the result differs from the standard basis $\{ \delta_w \mid w \in W\}$ of $H_A$ by a unitriangular matrix. This shows that $ch$ is an isomorphism of algebras. The inverse of $ch$ sends $b_s \mapsto [B_s]$. Finally, the double leaves basis for morphism spaces between Bott-Samelson objects is a homogeneous basis whose degrees match the defect formula by Lemma \ref{lem:degLL}, leading to the desired matching of sesquilinear forms.
\end{proof}


Applying Theorem \ref{thm:categorification theorem} to the '$p$-adapted realizations' from Corollary \ref{corollary: color graded Hecke category}, we obtain a categorification for a large class of unequal parameter Hecke algebras.

\begin{corollary}\label{corollary: categorification result for p grading}
Suppose that $\Bbbk$ has characteristic $p$ (possibly $0$) and let $V$ and $\HC$ be as in Corollary \ref{corollary: color graded Hecke category}. Then $[\HC]$ is isomorphic to the Hecke algebra of $W$ with generic unequal parameters $\{ v_s \mid s \in S\} $ subject to the relation $v_s = v_t$ whenever $s \sim_p t$.
\end{corollary}

\begin{remark} For any realization $V$ we have that $\sim_{\uneq}$ refines $\sim_{V}$, so there is a natural quotient homomorphism from $A_{\uneq}$ to $A_V := \Z^{S/\sim_V}$. We let
$H_V$ denote the corresponding specialization of $H_{\uneq}$. Recall from Example \ref{ex:universal} that there is always a realization $V'$ with the same Cartan matrix as $V$ (and
hence $\sim_V = \sim_{V'}$) for which $\HC(V')$ admits a grading by $\Z^{S/\sim_V}$. Using this we can categorify $H_V$. However, $\sim_V$ may be much coarser than $\sim_{\uneq}$, so
we may not be able to categorify $H_{\uneq}$. The previous corollary uses the fact that for some Coxeter groups, including finite and affine Weyl groups, there are realizations where $\sim_V$ agrees with $\sim_{\uneq}$.
\end{remark}

\subsection{Rescaled Hecke algebras} \label{subsec:rescaled}

So far we only considered gradings on the Hecke category for which $f_s = g_s$ for all $s \in S$. There are gradings which do not satisfy this property, but as we have seen in
Proposition \ref{prop:inner}, any grading on the Hecke category is inner relative to a grading which satisfies $f_s = g_s$. Thus, to determine the Grothendieck group of a general Hecke
category with a grading, we need to determine how the Grothendieck group changes along inner gradings. As one might expect, the Grothendieck group essentially does not change, but one
has to rescale objects. Note that both Assumption \ref{assume2} and Assumption \ref{assume3} can be dispensed with this way.

Suppose that $\HC$ is $A$-graded as in Assumptions \ref{assume} and \ref{assume2} and \ref{assume3}. Let $L$ be the corresponding $A$-valued parameter function. Let $d: A \rightarrow A'$ together with $\{ d(X) \in A' \mid X \in \HC \}$ be an inner grading datum and let $L' = d \circ L$. Because it is cumbersome to have parentheses within shifts, we write $d_X$ instead of $d(X)$, and $d_s$ instead of $d(B_s)$. As before we denote by $\HC_{A'}$ the category with objects $X(a')$ where $a' \in A'$ and $X\in \HC$ is a Bott-Samelson object.

As there is no definition of $H_{A'}$ in this generality, we instead use $\mathbb{Z}[A'] \otimes_{\mathbb{Z}[A]} H_A$ as a proxy. We write $\delta_s$ and $b_s$ for $1 \ot \delta_s$ and $1 \ot b_s$ respectively.


\begin{lemma}\label{lemma: inner grading on Grothendieck groups}
Under the assumptions above, there is an isomorphism of $\mathbb{Z}[A']$-algebras
\begin{align*}
	\mathbb{Z}[A'] \otimes_{\mathbb{Z}[A]} [\HC]_A &\overset{\sim}{\rightarrow} [\HC]_{A'} \\
	[X] & \mapsto [X(d_X)].
\end{align*}
In particular, there is an isomorphism  of $\mathbb{Z}[A']$-algebras
\begin{equation*}
\mathbb{Z}[A'] \otimes_{\mathbb{Z}[A]} H_A\overset{\sim}{\rightarrow} [\HC]_{A'}, \quad b_s \mapsto [B_s (d_s)].
\end{equation*}
\end{lemma}

\begin{proof}
We have a fully faithful functor $\HC_A \rightarrow \HC_{A'}$ which is given on objects by $X(a) \mapsto X(d(a) + d_X)$ and on morphisms by the canonical identification
\begin{align*}
\Hom_{\HC_A}(X(a), Y(b))& = \Hom^{b-a}_{\HC_A}(X,Y) \\
& \overset{def}{=} \Hom^{d(b-a) + d_Y-d_X}_{\HC_{A'}}(X,Y)\\
& = \Hom_{\HC_{L'}}(X(d(a) + d_X), Y(d(b) + d_Y)).
\end{align*}
Applying this map to the Grothendieck group yields a homomorphism of $\mathbb{Z}[A]$-algebras $ [\HC]_A \rightarrow [\HC]_{A'}$ where the $A$-action on the right hand side is via $d: A \rightarrow A'$. Hence, we get an induced algebra homomorphism $\mathbb{Z}[A'] \otimes_{\mathbb{Z}[A]} [\HC]_A \rightarrow [\HC]_{A'}$. This is an isomorphism since it identifies the basis $\{ [B_{\underline{w}} ] \mid w \in W \}$ with the rescaled basis $\{ v^{-d_{B_{\underline{w}}}} [B_{\underline{w}}] \mid w \in W \}$ (where we pick a reduced expression $\underline{w}$ for each $w \in W$).
\end{proof}

Under this isomorphism, the indecomposable object $B_s \in \HC_{A'}$ has image in $\mathbb{Z}[A'] \otimes_{\mathbb{Z}[A]} H_A$ equal to \begin{equation} b_s' := v^{d_s} b_s. \end{equation} A
common theme in Kazhdan-Lusztig theory is self-duality: for example, $b_s$ is fixed by the bar involution on $H_{A}$, and the object $B_s$ is self-dual under a natural duality
functor. In the usual $\Z$-grading, $B_s$ is the unique shift of itself which is self-dual. However, larger grading groups can have additional self-dual shifts. For example, the shift
$d_s = \frac{f_s - g_s}{2}$ used in Proposition \ref{prop:inner} is fixed under the bar involution, and consequently $b_s'$ is still self-dual! One can verify that whenever $d_s$ is
fixed by the bar involution, then the vertical flip of diagrams extends to an antiinvolution on $\HC$ which descends to the bar involution. If $\bar{d_s} \ne d_s$, then this statement
will fail. We leave the details to the interested reader.

Let $\{\delta_w\}$ be the standard basis of $H_A$, extended to $\mathbb{Z}[A'] \otimes_{\mathbb{Z}[A]} [\HC]_A$. Then $\HC_{A'}$ admits a rescaled variant on the character describing the isomorphism from Lemma \ref{lemma: inner grading on Grothendieck groups}:
\begin{equation}\label{eq: rescaled diagrammatic character}
\begin{aligned}
	{[B]} & \mapsto v^{d_B} \sum_{w} d(grk_R^A \Hom_{\HC_A, \ge w}(B_{\underline{w}}, B) )\delta_w \\
	& = \sum_{w} grk_R^{A'} \Hom_{\HC_{A'}, \ge w}(B_{\underline{w}}, B) v^{d_{B_{\underline{w}}}} \delta_w.
\end{aligned}
\end{equation}

To get a nicer formula for the character, it seems more natural to also replace the standard basis with a basis consisting of rescaled elements $\delta_{\underline{w}}' :=
v^{d_{B_{\underline{w}}}} \delta_w$, which might depend on the reduced expression $\underline{w}$. Let us briefly elaborate on what this rescaling does to the presentation of the Hecke
algebra.

The traditional definition of a generic Hecke algebra (see e.g. \cite[\S 7]{HumpCox}) is an algebra with a specific kind of presentation. Namely, it is generated by elements $\delta_s$ for $s \in S$, subject to the braid relations and to the quadratic relation
\begin{equation} \label{genericquadratic} \delta_s^2 = p_s \delta_s + q_s \end{equation}
for elements $p_s, q_s$ in the base ring which are constant on equivalence classes for $\sim_{\uneq}$.

Now let $\delta_s' := v_s \delta_s$ for some invertible parameter $v_s$. The quadratic relation
satisfied by $\delta_s'$ is now \begin{equation} (\delta_s')^2 = v_s p_s \delta_s' + v_s^2 q_s. \end{equation} Thus $\{\delta_s'\}$ also satisfy the generic quadratic relation but with different parameters $p_s' := v_s p_s$ and $q_s' := v_s^2 q_s$. If $v_s = v_t$ whenever $s \sim_{\uneq} t$, then $p_s'$ and $q_s'$ are constant on equivalence classes, and $\{\delta_s'\}$ continue to satisfy the braid relation. This gives two Hecke
presentations for the same Hecke algebra, one with generators $\delta_s$ and the other with generators $\delta'_s$. Thus under Assumption \ref{assume2} we can directly describe $[\HC]_{A'}$ as a generic Hecke algebra, already encompassing its rescaling into the notation.

However, if $v_s \ne v_t$ but $m_{st}$ is odd (e.g. if Assumption \ref{assume2} does not hold), then $\delta_s'$ and $\delta_t'$ do not satisfy the braid relation. Instead, one has a more general presentation with relation
\begin{equation*}
	\delta'_s\delta'_t\delta'_s ... = v^{d_s - d_t} \delta'_t\delta'_s\delta'_t...
\end{equation*}
Trying to impose the usual braid relation when $p_s \ne p_t$ or $q_s \ne q_t$ for $m_{st}$ odd will result in an algebra with unusual relations, which may not have the same size as $W$.

The literature on algebras with these variations on the generic Hecke presentation is not developed, largely since everything can be obtained by isomorphism with a generic Hecke
algebra. The reader interested in abandoning Assumption \ref{assume2} should be prepared to deal with these technicalities, though the cost of abandoning Assumption
\ref{assume3} is not so high.

\begin{remark} \label{rmk:degrees of light leaves} When Assumption \ref{assume2} does not hold, one can still compute the character using \eqref{eq: rescaled diagrammatic character},
or compute the degree of a light leaf map, though it includes an extra scaling adjustment depending on the various reduced expressions involved. This adjustment is similar to the
calculation from Lemma \ref{lemma:secret inner} using the integers $n_s$, and we leave it to the reader. \end{remark}

\section{An example: The infinite dihedral group}\label{section:dihedral example}

Throughout this section assume that $\Bbbk$ is a field of characteristic $0$. Let
\begin{equation*}
    W = D_{\infty} = \langle s_1 , s_2 \rangle
\end{equation*}
be the infinite dihedral group. Let $\HC$ be the Hecke category associated to the realization
\begin{equation*}
    V:= \Bbbk \alpha_1 \oplus \Bbbk \alpha_2
\end{equation*}
where $\langle \alpha_i, \alpha_j^{\vee} \rangle = 2 \delta_{ij}$ for $i,j \in \{s,t\}$. This is the inflation to the infinite dihedral group of the reflection representation in type $A_1 \times A_1$. We equip $V$ with a grading by $\Z^2$ where $\deg(\alpha_1) = (2,0)$ and $\deg(\alpha_2) = (0,2)$, and we equip $\HC$ with a grading so that $s_1$-colored dots have degree $(1,0)$ and $s_2$-colored dots have degree $(0,1)$. In particular, $f_s = g_s$ for all $s \in S$, so Assumption \ref{assume3} holds.

By Theorem \ref{thm:categorification theorem} the Grothendieck group $H$ of $\HC$ is the Hecke algebra with unequal parameters, over the base ring $\Z[v_1, v_1^{-1}, v_2, v_2^{-1}]$. We want to compute the basis for the Hecke algebra coming from indecomposable objects $\HC$, which we call the \emph{double-0 canonical basis}. To do so we need to better understand Jones-Wenzl projectors in $\TL_n(0)$.


\subsection{Two-step recursion for Jones-Wenzl projectors}

\begin{lemma}
    Assume the base ring $\Bbbk$ has characteristic zero. The Jones-Wenzl projector exists in $\TL_n(0)$ if and only if $n$ is odd.
\end{lemma}
\begin{proof}
    It is well-known that the Jones-Wenzl projector $\JW_n$ exists if and only ${n\brack k}$ is invertible in $\Bbbk$ for all $0 \le k \le n$ (c.f. \cite[Theorem A]{HaziRotatable}). By \eqref{quantumbinomclaim} this is equivalent to $n$ being odd.
\end{proof}

When $\JW_n$ exists one can define the scalars $p_1^{(n)}$ and $p_2^{(n)}$ in $\Bbbk$, as coefficients of the identity in one-step or two-step partial traces, see \cite[(6.40), (6.41)]{EWLoc}. When both $\JW_{n+2}$ and $\JW_n$ exist, one has
\begin{equation}
    \vcenter{\hbox{\begin{tikzpicture}[baseline={(current bounding box.center)}, scale=0.6]
    
    \draw (0,0) rectangle (2,1) node[pos=.5] {$\JW_{n+2}$};
    
    \draw (1.8,1) to[out=90,in=180] (2,1.2) to[out=0,in=90] (2.2,1) to (2.2,0) to[out=270,in=0] (2,-0.2) to[out=180,in=270] (1.8,0);
    \draw (1.6,1) to[out=90,in=180] (2,1.4) to[out=0,in=90] (2.4,1) to (2.4,0) to[out=270,in=0] (2,-0.4) to[out=180,in=270] (1.6,0);
    
    \draw (0.2,1) to (0.2, 1.4);
    \draw (1.4,1) to (1.4, 1.4);
    \draw (0.2,0) to (0.2, -0.4);
    \draw (1.4,0) to (1.4, -0.4);
    
    \node at (0.5,1.2) [circle,fill,inner sep=.5pt]{};
    \node at (0.8,1.2) [circle,fill,inner sep=.5pt]{};
    \node at (1.1,1.2) [circle,fill,inner sep=.5pt]{};
    \node at (0.5,-0.2) [circle,fill,inner sep=.5pt]{};
    \node at (0.8,-0.2) [circle,fill,inner sep=.5pt]{};
    \node at (1.1,-0.2) [circle,fill,inner sep=.5pt]{};
    \end{tikzpicture}}}
    = p_2^{(n+2)} \cdot
    \vcenter{\hbox{\begin{tikzpicture}[baseline={(current bounding box.center)}, scale=0.6]
    
    \draw (0,0) rectangle (2,1) node[pos=.5] {$\JW_{n}$};

    \draw (0.2,1) to (0.2, 1.4);
    \draw (1.8,1) to (1.8, 1.4);
    \draw (0.2,0) to (0.2, -0.4);
    \draw (1.8,0) to (1.8, -0.4);
    
    \node at (0.6,1.2) [circle,fill,inner sep=.5pt]{};
    \node at (1.0,1.2) [circle,fill,inner sep=.5pt]{};
    \node at (1.4,1.2) [circle,fill,inner sep=.5pt]{};
    \node at (0.6,-0.2) [circle,fill,inner sep=.5pt]{};
    \node at (1.0,-0.2) [circle,fill,inner sep=.5pt]{};
    \node at (1.4,-0.2) [circle,fill,inner sep=.5pt]{};
    \end{tikzpicture}}}.
\end{equation}
It is not obvious that $p_2^{(n+2)}$ should be invertible, but we will prove below that it is.

In \cite[Theorem 6.33]{EWLoc} a two-step recursion for Jones-Wenzl projectors is developed, where one can compute $\JW_{n+2}$ using $\JW_n$ without assuming that $\JW_{n+1}$ is well-defined. More precisely, when $\JW_n$ exists, \cite[Theorem 6.33]{EWLoc} defines a morphism $J \in \TL_{n+2}$ and gives a precise criterion for when $J = \JW_{n+2}$. Note that \cite[Theorem 6.33]{EWLoc} is not correctly stated for $n<2$, as $p_2^{(n)}$ is not well-defined (for $n=1$ the terms and criteria involving $p_2^{(n)}$ or its inverse should simply be omitted).

The formulas for $J$ and the criteria simplify dramatically under the hypothesis that $[2]_{s,t} = [2]_{t,s} = 0$. Rather than reproduce \cite[Theorem 6.33]{EWLoc} here in all its glory, we examine the consequences. For $n \ge 1$ odd let
\begin{equation}
    E_0 = E_0^{(n)} :=
    \vcenter{\hbox{\begin{tikzpicture}[baseline={(current bounding box.center)}, scale=0.6]
    
    \draw (0,0) rectangle (2,1) node[pos=.5] {$\JW_{n}$};

    \draw (0.2,1) to (0.2, 1.4);
    \draw (1.8,1) to (1.8, 1.4);
    \draw (0.2,0) to (0.2, -0.4);
    \draw (1.8,0) to (1.8, -0.4);
    
    \node at (0.6,1.2) [circle,fill,inner sep=.5pt]{};
    \node at (1.0,1.2) [circle,fill,inner sep=.5pt]{};
    \node at (1.4,1.2) [circle,fill,inner sep=.5pt]{};
    \node at (0.6,-0.2) [circle,fill,inner sep=.5pt]{};
    \node at (1.0,-0.2) [circle,fill,inner sep=.5pt]{};
    \node at (1.4,-0.2) [circle,fill,inner sep=.5pt]{};

    \draw (2.3,1.4) to (2.3, -0.4);
    \draw (2.6,1.4) to (2.6, -0.4);
    \end{tikzpicture}}},
    %
    %
    %
    \quad
    E_1 = E_1^{(n)} :=
    \vcenter{\hbox{\begin{tikzpicture}[baseline={(current bounding box.center)}, scale=0.6]
    
    \draw (0,0) rectangle (2,1) node[pos=.5] {$\JW_{n}$};

    \draw (0,2) rectangle (2,3) node[pos=.5] {$\JW_{n}$};

    \draw (0.2,1) to (0.2, 2);
    \draw (1.4,1) to (1.4, 2);
    \draw (0.2,0) to (0.2, -0.4);
    \draw (1.8,0) to (1.8, -0.4);
    \draw (0.2,3) to (0.2, 3.4);
    \draw (1.8,3) to (1.8, 3.4);

    \node at (0.6,3.2) [circle,fill,inner sep=.5pt]{};
    \node at (1.0,3.2) [circle,fill,inner sep=.5pt]{};
    \node at (1.4,3.2) [circle,fill,inner sep=.5pt]{};
    \node at (0.5,1.5) [circle,fill,inner sep=.5pt]{};
    \node at (0.8,1.5) [circle,fill,inner sep=.5pt]{};
    \node at (1.1,1.5) [circle,fill,inner sep=.5pt]{};
    \node at (0.6,-0.2) [circle,fill,inner sep=.5pt]{};
    \node at (1.0,-0.2) [circle,fill,inner sep=.5pt]{};
    \node at (1.4,-0.2) [circle,fill,inner sep=.5pt]{};

    \draw (2.3,-0.4) to[out=90,in=180] (2.45,0) to[out=0,in=90] (2.6,-0.4);
    \draw (1.8, 1) to[out=90, in=180] (2.2, 1.5) to[out=0,in=270] (2.6,2) to (2.6,3.4);
    \draw (1.8,2) to[out=270,in=180] (2.05,1.8) to[out=0,in=270] (2.3,2) to (2.3,3.4);

    \end{tikzpicture}}},
    \quad
    E_2 = E_2^{(n)} := \vcenter{\hbox{\begin{tikzpicture}[baseline={(current bounding box.center)}, scale=0.6]
    
    \draw (0,0) rectangle (2,1) node[pos=.5] {$\JW_{n}$};

    \draw (0,2) rectangle (2,3) node[pos=.5] {$\JW_{n}$};

    \draw (0.2,1) to (0.2, 2);
    \draw (1.4,1) to (1.4, 2);
    \draw (0.2,0) to (0.2, -0.4);
    \draw (1.8,0) to (1.8, -0.4);
    \draw (0.2,3) to (0.2, 3.4);
    \draw (1.8,3) to (1.8, 3.4);

    \node at (0.6,3.2) [circle,fill,inner sep=.5pt]{};
    \node at (1.0,3.2) [circle,fill,inner sep=.5pt]{};
    \node at (1.4,3.2) [circle,fill,inner sep=.5pt]{};
    \node at (0.5,1.5) [circle,fill,inner sep=.5pt]{};
    \node at (0.8,1.5) [circle,fill,inner sep=.5pt]{};
    \node at (1.1,1.5) [circle,fill,inner sep=.5pt]{};
    \node at (0.6,-0.2) [circle,fill,inner sep=.5pt]{};
    \node at (1.0,-0.2) [circle,fill,inner sep=.5pt]{};
    \node at (1.4,-0.2) [circle,fill,inner sep=.5pt]{};

    \draw (2.3,3.4) to[out=90,in=180] (2.45,3) to[out=0,in=90] (2.6,3.4);
    \draw (1.8, 2) to[out=270, in=180] (2.2, 1.5) to[out=0,in=90] (2.6,1) to (2.6,-0.4);
    \draw (1.8,1) to[out=90,in=180] (2.05,1.2) to[out=0,in=90] (2.3,1) to (2.3,-0.4);

    \end{tikzpicture}}},
\end{equation}
and for $n \ge 3$ odd let
\begin{equation}
    E_3 = E_3^{(n)} :=  \text{ } \vcenter{\hbox{\begin{tikzpicture}[baseline={(current bounding box.center)}, scale=0.6]
    
    \draw (0,0) rectangle (2,1) node[pos=.5] {$\JW_{n}$};

    \draw (0,2) rectangle (2,3) node[pos=.5] {$\JW_{n}$};

    \draw (0.2,1) to (0.2, 2);
    \draw (1.4,1) to (1.4, 2);
    \draw (0.2,0) to (0.2, -0.4);
    \draw (1.8,0) to (1.8, -0.4);
    \draw (0.2,3) to (0.2, 3.4);
    \draw (1.8,3) to (1.8, 3.4);

    \node at (0.6,3.2) [circle,fill,inner sep=.5pt]{};
    \node at (1.0,3.2) [circle,fill,inner sep=.5pt]{};
    \node at (1.4,3.2) [circle,fill,inner sep=.5pt]{};
    \node at (0.5,1.5) [circle,fill,inner sep=.5pt]{};
    \node at (0.8,1.5) [circle,fill,inner sep=.5pt]{};
    \node at (1.1,1.5) [circle,fill,inner sep=.5pt]{};
    \node at (0.6,-0.2) [circle,fill,inner sep=.5pt]{};
    \node at (1.0,-0.2) [circle,fill,inner sep=.5pt]{};
    \node at (1.4,-0.2) [circle,fill,inner sep=.5pt]{};

    \draw (1.6,2) to[out=270,in=180] (2,1.6) to[out=0,in=270] (2.4,2) to (2.4,3.4);
    \draw (1.8,2) to[out=270,in=180] (2,1.8) to[out=0,in=270] (2.2,2) to (2.2,3.4);
    \draw (1.6,1) to[out=90,in=180] (2,1.4) to[out=0,in=90] (2.4,1) to (2.4,-0.4);
    \draw (1.8,1) to[out=90,in=180] (2,1.2) to[out=0,in=90] (2.2,1) to (2.2,-0.4);
    \end{tikzpicture}}} \; .
    \end{equation}

\begin{thm} Assume the base ring $\Bbbk$ has characteristic zero. Let $n \ge 1$ be odd. Then $p_1^{(n)} = 0$. If $n \ge 3$ then $p_2^{(n)}$ is invertible, and satisfies the recursive formula
\begin{equation} \label{eq:p2recursive} p_2^{(3)} = -2, \qquad p_2^{(n+2)} = -2 - \frac{1}{p_2^{(n)}}. \end{equation}
and the closed formula
\begin{equation} \label{eq:p2closed} p_2^{(2k+1)} = \frac{-(k+1)}{k}. \end{equation}
We have
\begin{equation} \JW_3 = E_0^{(1)} - E_1^{(1)} - E_2^{(1)}. \end{equation}
When $n \ge 3$ we have
\begin{equation} \label{eq:JWrecursive2} \JW_{n+2} = E_0^{(n)} - E_1^{(n)} - E_2^{(n)} - \frac{1}{p_2^{(n)}} E_3^{(n)}. \end{equation}
Moreover, the elements $\JW_{n+2}$, $E_1^{(n)}$, $E_2^{(n)}$, and $\frac{1}{p_2^{(n)}} E_3^{(n)}$ are pairwise orthogonal idempotents in $\TL_{n+2}(0)$, giving a decomposition of the idempotent $E_0^{(n)}$, as in
\begin{equation} \label{eq:E0decomp} E_0^{(n)} = \JW_{n+2} + E_1^{(n)} + E_2^{(n)} + \frac{1}{p_2^{(n)}} E_3^{(n)}. \end{equation}
\end{thm}

\begin{proof}
These statements can be verified for small $n$ directly, so we focus on the inductive step, using \cite[Theorem 6.33]{EWLoc}. We assume that $p_1^{(n)} = 0$ and $p_2^{(n)}$ obeys \eqref{eq:p2closed}. Then $p_2^{(n)}$ and $1 + [2]_{s,t} p_1^{(n)}$ are both invertible. So the assumptions of the second half of \cite[Theorem 6.33]{EWLoc} are active, and $\kappa = -1$ since $[2]_{s,t} = 0$. By \cite[(6.45)]{EWLoc} we have $p_1^{(n+2)} = - \kappa p_1^{(n)} = 0$ and 
\[ p_2^{(n+2)} = (1 - \kappa)[3] + \frac{\kappa}{p_2^{(n)}} = -2-\frac{1}{p_2^{(n)}},\]
since $[3] = -1$. This proves the inductive step for \eqref{eq:p2recursive}, from which \eqref{eq:p2closed} follows easily. The formula for $\JW_{n+2}$ given in the second half of \cite[Theorem 6.33]{EWLoc} easily specializes to the formula \eqref{eq:JWrecursive2} when $\kappa = -1$ and $[2]_{s,t} = p_1^{(n)} = 0$. The statements about orthogonal idempotents are easy to verify, and \eqref{eq:E0decomp} is just a rearrangement of \eqref{eq:JWrecursive2}. \end{proof}

\begin{rem}
    It should be possible to also deduce the recursion \eqref{eq:JWrecursive2} from the results in \cite[Thm 4.8]{SuttonTiltMod}. Note however that the statement of \cite[Thm 4.8]{SuttonTiltMod} is not quite correct as stated since the second term in the definition of $A_v^i$ (\cite[Definition 4.7]{SuttonTiltMod}) is missing a coefficient. We thank Daniel Tubbenhauer and Paul Wedrich for pointing this out to us. For $\delta = 0$ these coefficients can be chosen to be $0$ and applying the recursion from \cite[Thm 4.8]{SuttonTiltMod} twice results in the recursion from \eqref{eq:JWrecursive2} above.
\end{rem}

The theorem above can be used to give a proof of the homogeneity of all Jones-Wenzl projectors (again we emphasize that we work in characteristic $0$ here so Proposition \ref{prop: Homogenity of Jones-Wenzl} does not apply).
\begin{corollary}\label{cor: homogeneity of JW in characteristic 0}
    In $\TL_n(0)$ the Jones-Wenzl projectors $\JW_n$ for $n$ odd are all homogeneous of degree $0$. That is, the number of even caps equals the number of odd cups.
\end{corollary}
\begin{proof}
    For $n=1,3$ this can be checked by direct inspection:
    \begin{equation*}
        \JW_1 = 
        \vcenter{\hbox{\begin{tikzpicture}[baseline={(current bounding box.center)}, scale=0.6]
            \draw (0,0) to (0,1);
        \end{tikzpicture}}} \text{ },
        \quad
        \JW_3 = 
        \vcenter{\hbox{\begin{tikzpicture}[baseline={(current bounding box.center)}, scale=0.6]
            \draw (0,0) to (0,1);
            \draw (0.3,0) to (0.3,1);
            \draw (0.6,0) to (0.6,1);
        \end{tikzpicture}}}
        \text{ } - \text{ }
        \vcenter{\hbox{\begin{tikzpicture}[baseline={(current bounding box.center)}, scale=0.6]
            \draw (0,0) to[out=90,in=180] (0.15,0.2) to[out=0,in=90] (0.3,0);
            \draw (0,1) to[out=270,in=180] (0.3,0.5) to[out=0,in=90] (0.6,0);
            \draw (0.3,1) to[out=270,in=180] (0.45,0.8) to[out=0,in=270](0.6,1);
        \end{tikzpicture}}}
        \text{ } - \text{ }
        \vcenter{\hbox{\begin{tikzpicture}[baseline={(current bounding box.center)}, scale=0.6]
            \draw (0.3,0) to[out=90,in=180] (0.45,0.2) to[out=0,in=90] (0.6,0);
            \draw (0,0) to[out=90,in=180] (0.3,0.5) to[out=0,in=270] (0.6,1);
            \draw (0,1) to[out=270,in=180] (0.15,0.8) to[out=0,in=270](0.3,1);
        \end{tikzpicture}}}.
    \end{equation*}
    For $n >3$ this follows by induction using \eqref{eq:JWrecursive2}.
\end{proof}

Finally we note the following indecomposability result, which is also an easy corollary of the fact that even Jones-Wenzl projectors do not exist.

\begin{corollary} \label{cor:indecompy} Let $n \ge 1$ be odd. In $\TL_{n+1}(0)$ the idempotent
	$y_{n+1} := \vcenter{\hbox{\begin{tikzpicture}[baseline={(current bounding box.center)}, scale=0.6]
    
    \draw (0,0) rectangle (2,1) node[pos=.5] {$\JW_{n}$};

    \draw (0.2,1) to (0.2, 1.4);
    \draw (1.8,1) to (1.8, 1.4);
    \draw (0.2,0) to (0.2, -0.4);
    \draw (1.8,0) to (1.8, -0.4);
    
    \node at (0.6,1.2) [circle,fill,inner sep=.5pt]{};
    \node at (1.0,1.2) [circle,fill,inner sep=.5pt]{};
    \node at (1.4,1.2) [circle,fill,inner sep=.5pt]{};
    \node at (0.6,-0.2) [circle,fill,inner sep=.5pt]{};
    \node at (1.0,-0.2) [circle,fill,inner sep=.5pt]{};
    \node at (1.4,-0.2) [circle,fill,inner sep=.5pt]{};

    \draw (2.3,1.4) to (2.3, -0.4);
    \end{tikzpicture}}}$
	is primitive. \end{corollary}

\begin{proof} By the usual arguments one has
    \begin{equation*}
        y_{n+1} \TL_{n+1}(0) y_{n+1} = \text{Span}_{\Bbbk} \left\{ \; \vcenter{\hbox{\begin{tikzpicture}[baseline={(current bounding box.center)}, scale=0.6]
    
    \draw (0,0) rectangle (2,1) node[pos=.5] {$\JW_{n}$};

    \draw (0.2,1) to (0.2, 1.4);
    \draw (1.8,1) to (1.8, 1.4);
    \draw (0.2,0) to (0.2, -0.4);
    \draw (1.8,0) to (1.8, -0.4);
    
    \node at (0.6,1.2) [circle,fill,inner sep=.5pt]{};
    \node at (1.0,1.2) [circle,fill,inner sep=.5pt]{};
    \node at (1.4,1.2) [circle,fill,inner sep=.5pt]{};
    \node at (0.6,-0.2) [circle,fill,inner sep=.5pt]{};
    \node at (1.0,-0.2) [circle,fill,inner sep=.5pt]{};
    \node at (1.4,-0.2) [circle,fill,inner sep=.5pt]{};

    \draw (2.3,1.4) to (2.3, -0.4);
    \end{tikzpicture}}} \;, \; \;
        \vcenter{\hbox{\begin{tikzpicture}[baseline={(current bounding box.center)}, scale=0.6]
    
    \draw (0,0) rectangle (2,1) node[pos=.5] {$\JW_{n}$};

    \draw (0,2) rectangle (2,3) node[pos=.5] {$\JW_{n}$};

    \draw (0.2,1) to (0.2, 2);
    \draw (1.4,1) to (1.4, 2);
    \draw (0.2,0) to (0.2, -0.4);
    \draw (1.8,0) to (1.8, -0.4);
    \draw (0.2,3) to (0.2, 3.4);
    \draw (1.8,3) to (1.8, 3.4);

    \node at (0.6,3.2) [circle,fill,inner sep=.5pt]{};
    \node at (1.0,3.2) [circle,fill,inner sep=.5pt]{};
    \node at (1.4,3.2) [circle,fill,inner sep=.5pt]{};
    \node at (0.5,1.5) [circle,fill,inner sep=.5pt]{};
    \node at (0.8,1.5) [circle,fill,inner sep=.5pt]{};
    \node at (1.1,1.5) [circle,fill,inner sep=.5pt]{};
    \node at (0.6,-0.2) [circle,fill,inner sep=.5pt]{};
    \node at (1.0,-0.2) [circle,fill,inner sep=.5pt]{};
    \node at (1.4,-0.2) [circle,fill,inner sep=.5pt]{};

    \draw (1.8,1) to[out=90,in=180] (2.05,1.2) to[out=0,in=90] (2.3,1) to (2.3,-0.4);
    \draw (1.8,2) to[out=270,in=180] (2.05,1.8) to[out=0,in=270] (2.3,2) to (2.3,3.4);

    \end{tikzpicture}}} \;
        \right\}.
    \end{equation*}
	Using that $p_1^{(n)} = 0$, this is isomorphic to $\Bbbk[x]/(x^2)$ as an algebra, with identity element $y_{n+1}$. Thus $y_{n+1}$ is a primitive idempotent.
\end{proof}

\subsection{The basis from indecomposable objects}

For each $k \ge 1$ let
\begin{align*}
    1_k &:= s_1 s_2s_1... \\
    2_k &:= s_2 s_1 s_2 ...
\end{align*}
be the unique elements of length $k$ starting in $s_1$ and $s_2$ respectively. Let $b_{1_k}$ and $b_{2_k}$ denote the images in the Hecke algebra of the symbols of indecomposable objects in $\HC$, also called the double-0 canonical basis.

\begin{lemma}\label{lemma: recursive formula for p-can base of even length}
    We have
    \begin{align}
        &b_{1}b_{2} = b_{1_2},  & & b_{2}b_{1} = b_{2_2},\\
        &b_{1}b_{2} b_1 = b_{1_3}, \qquad &  &b_{2} b_{1} b_2 = b_{2_3}, \nonumber \\
        &b_{1}b_{2} b_{1} b_{2}  =b_{1_4} + (v_1 v_2^{-1} + v_1^{-1} v_2) b_{1_2}, \qquad & & b_2 b_1 b_2 b_1 =  b_{2_4} + (v_1^{-1} v_2  + v_1 v_2^{-1})b_{2_2}.   \nonumber
    \end{align}
    Moreover, for any $k \ge 2$ we have
    \begin{align}
		& b_{1_{2k}}b_1 = b_{1_{2k+1}}, \qquad b_{2_{2k}}b_2 = b_{2_{2k+1}}, \\
        & b_{1_{2k}}b_{1}b_{2} = b_{1_{2k+2}}+ (v_1 v_2^{-1} + v_1^{-1}v_2 )  b_{1_{2k}} +  b_{1_{2k-2}}, \nonumber \\
        & b_{2_{2k}}b_{2}b_{1} =  b_{2_{2k+2}}+ (v_1 v_2^{-1} + v_1^{-1}v_2 )  b_{2_{2k}} +  b_{2_{2k-2}}. \nonumber
    \end{align}
\end{lemma}

\begin{proof}
    We have $\End(1_2) = \End(2_2) = S(V)$ which shows that $1_2$ and $2_2$ are indecomposable. Thus $b_{1_2} = b_1 b_2$ and $b_{2_2} = b_2 b_1$.
    
    Any degree $0$ morphism in $\End(1_4)$ in the standard $\mathbb{Z}$-grading comes from the Temperley-Lieb algebra \cite[Prop. 1.1]{ECathedral}. Since our $\mathbb{Z}^2$-grading refines this grading, any degree $0$ morphism in $\End(1_4)$ in the $\mathbb{Z}^2$-grading also comes from the Temperley-Lieb algebra, though not all Temperley-Lieb diagrams have degree $0$. In $TL_3(0)$ we have the decomposition into primitive idempotents
    \begin{equation}\label{eq: idempotent dec for JW3}
        \vcenter{\hbox{\begin{tikzpicture}[baseline={(current bounding box.center)}, scale=0.6]
            \draw (0,0) to (0,1);
            \draw (0.3,0) to (0.3,1);
            \draw (0.6,0) to (0.6,1);
        \end{tikzpicture}}}
        \text{ } = \text{ } \JW_3 + \text{ }
        \vcenter{\hbox{\begin{tikzpicture}[baseline={(current bounding box.center)}, scale=0.6]
            \draw (0,0) to[out=90,in=180] (0.15,0.2) to[out=0,in=90] (0.3,0);
            \draw (0,1) to[out=270,in=180] (0.3,0.5) to[out=0,in=90] (0.6,0);
            \draw (0.3,1) to[out=270,in=180] (0.45,0.8) to[out=0,in=270](0.6,1);
        \end{tikzpicture}}}
        \text{ } + \text{ }
        \vcenter{\hbox{\begin{tikzpicture}[baseline={(current bounding box.center)}, scale=0.6]
            \draw (0.3,0) to[out=90,in=180] (0.45,0.2) to[out=0,in=90] (0.6,0);
            \draw (0,0) to[out=90,in=180] (0.3,0.5) to[out=0,in=270] (0.6,1);
            \draw (0,1) to[out=270,in=180] (0.15,0.8) to[out=0,in=270](0.3,1);
    \end{tikzpicture}}}.
    \end{equation}
    Since all terms in \eqref{eq: idempotent dec for JW3} are homogeneous of degree $0$, it follows that \eqref{eq: idempotent dec for JW3} induces a decomposition of $1$ into primitive idempotents in $\End^0(1_4)$ (with respect to the $\mathbb{Z}^2$-grading). It is straightforward to check that the idempotent corresponding to $\vcenter{\hbox{\begin{tikzpicture}[baseline={(current bounding box.center)}, scale=0.6]
            \draw (0,0) to[out=90,in=180] (0.15,0.2) to[out=0,in=90] (0.3,0);
            \draw (0,1) to[out=270,in=180] (0.3,0.5) to[out=0,in=90] (0.6,0);
            \draw (0.3,1) to[out=270,in=180] (0.45,0.8) to[out=0,in=270](0.6,1);
        \end{tikzpicture}}}$
    picks out the object $B_{1_2} (g_2 - f_1)$ in the Hecke category and the idempotent corresponding to
    $\vcenter{\hbox{\begin{tikzpicture}[baseline={(current bounding box.center)}, scale=0.6]
            \draw (0.3,0) to[out=90,in=180] (0.45,0.2) to[out=0,in=90] (0.6,0);
            \draw (0,0) to[out=90,in=180] (0.3,0.5) to[out=0,in=270] (0.6,1);
            \draw (0,1) to[out=270,in=180] (0.15,0.8) to[out=0,in=270](0.3,1);
    \end{tikzpicture}}}$
    picks out $B_{1_2}(g_1 - f_2)$. It follows from this that in the Hecke category we have
    \begin{equation*}
        1_4 \cong X \oplus B_{1_2}(-1,1) \oplus B_{1_2}(1,-1)
    \end{equation*}
    where $X $ is the indecomposable summand corresponding to $\JW_3$. Hence, we must have $X \cong B_{1_4}$ and taking the Grothendieck group yields
    \begin{equation*}
         b_1 b_2 b_1 b_2  = = b_{1_4} + (v_1 v_2^{-1} + v_1^{-1} v_2 )b _1 b_2.
     \end{equation*}
    Swapping the roles of $s_1$ and $s_2$ yields $b_2 b_1 b_2 b_1 =  b_{2_4} + (v_1 v_2^{-1} + v_1^{-1} v_2)b_2b_1$.

    Now let $k \ge 2 $. As before any degree $0$ morphism (with respect to the $\mathbb{Z}^2$-grading) in $\End(1_{2k})$ comes from the Temperley-Lieb algebra. Using \eqref{eq:JWrecursive2} with $n = 2k-1$ we see that
    \begin{equation*}
        Y B_1 B_2 = X_0 \oplus X_1  \oplus X_2   \oplus X_4.
    \end{equation*}
    where $X_i$ is the object corresponding to $E_i$ and $Y$ is the object corresponding to $\JW_{2k-1}$. Note that all the $E_i$ are homogeneous of degree $0$ by Corollary \ref{cor: homogeneity of JW in characteristic 0}. Clearly $X_0$ corresponds to the indecomposable object defined by $\JW_{2k+1}$. Moreover, one can check that $X_1 \cong Y(1,-1)$ and $X_2 \cong Y(-1,1)$. Finally, $X_3$ is the indecomposable object corresponding to $\JW_{2k-3}$. One can deduce from this by induction that the $\JW_{2k-1}$ pick out the objects $B_{2k}$ for all $k$ and for $k \ge 2$ we have
    \begin{equation*}
        B_{1_{2k}} B_1 B_2 \cong B_{1_{2k+2}} \oplus  B_{1_{2k}}(1,-1) \oplus   B_{1_{2k}}(-1,1) \oplus B_{1_{2k-2}}
    \end{equation*}
and thus $b_{1_{2k}}b_{1}b_{2} = b_{1_{2k+2}}+ (v_1 v_2^{-1} + v_1^{-1} v_2)  b_{1_{2k}} +  b_{1_{2k-2}}$ by taking the Grothendieck group.

Similarly, Corollary \ref{cor:indecompy} implies that $Y B_1$ is indecomposable, being the image of the primitive idempotent $y_{2k}$. This yields the remaining formulas.
\end{proof}

Let us describe the double-0 canonical basis in terms of the standard basis.

\begin{defn} For $w = s_1 s_2 s_1 ...$ (resp. $w = s_2 s_1 s_2 ...$) we write
\begin{equation}
	l_w =(1,0) + (0,1) + (1,0) + ... \quad (\text{resp. } l_w = (0,1) + (1,0) + (0,1) + ...).
\end{equation}
Let $\Gamma_w := \sum_{y \le w} v^{l_w - l_y} \delta_y$ in $H$. \end{defn}

\begin{corollary}
    For any $k \ge 1$ we have
    \begin{align}
        & b_{1_{2k}} = \Gamma_{1_{2k}}, \qquad & & b_{2_{2k}} = \Gamma_{2_{2k}}, \\
        & b_{1_{2k+1}} = \Gamma_{1_{2k+1}} + v_1^{-1} v_2 \Gamma_{1_{2k-1}}, \qquad & &
        b_{2_{2k+1}} = \Gamma_{2_{2k+1}} + v_1 v_2^{-1} \Gamma_{2_{2k-1}}. \nonumber
    \end{align}
\end{corollary}

\begin{proof} We leave the base cases ($k = 1$) to the reader.
    A straightforward computation shows that for $k \ge 1$ we have
	\begin{align*}
		\Gamma_{1_{2k}} b_1 &= \Gamma_{1_{2k+1}} + v_1^{-1} v_2 \Gamma_{1_{2k-1}} \\
		\Gamma_{1_{2k+1}} b_2 &= \Gamma_{1_{2k+2}} + v_1 v_2^{-1} \Gamma_{1_{2k}}.
	\end{align*}	    
    Thus, for $k \ge 2$ we have have
    \begin{align*}
        \Gamma_{1_{2k}} b_1 b_2 &= (\Gamma_{1_{2k+1}} + v_1^{-1} v_2 \Gamma_{1_{2k-1}}) b_2 \\
        &=   \Gamma_{2k+2} + v_1 v_2^{-1}  \Gamma_{2k} + v_1^{-1} v_2 (   \Gamma_{2k} + v_1 v_2^{-1}  \Gamma_{2k-2} )\\
        &= \Gamma_{2k+2} + (v_1 v_2^{-1} + v_1 ^{-1} v_2 ) \Gamma_{2k} + \Gamma_{2k-2}.
    \end{align*}
    It now follows from Lemma \ref{lemma: recursive formula for p-can base of even length} by induction that
    \begin{equation*}
        b_{1_{2k}} = \Gamma_{1_{2k}}
    \end{equation*}
    for all $k \ge 1$. Finally, we have
    \begin{equation*}
        b_{1_{2k+1}} = b_{1_{2k}} b_1 = \Gamma_{1_{2k}} b_1 =  \Gamma_{1_{2k+1}} + v_1^{-1} v_2 \Gamma_{1_{2k-1}}.
    \end{equation*}
    Swapping the role of $s_1$ and $s_2$ yields the corresponding formulas for $b_{2_{2k}}$ and $b_{2_{2k+1}}$.
\end{proof}

The computation above has the following interesting consequence in cell theory.

\begin{cor} The Hecke algebra with unequal parameters and the basis given above has 5 right cells: $\{\id\}$, $\{s_1\}$, $\{s_2\}$, $\{1_k\}_{k \ge 2}$, and $\{2_k\}_{k \ge 2}$. It has 4 two-sided cells: the last two right cells are merged into a two-sided cell, and the others remain as singletons. \end{cor}

Thus the two additional one-dimensional representations of the Hecke algebra (other than the sign and trivial representations) are actually cell representations.

For the remarks below we use the following notation. In the Karoubi envelope of the category $\TL(0)$, let $T_n$ denote the top summand of the object $n$. It is defined as the image of the unique primitive idempotent in $\End(n)$ containing the identity diagram with nonzero coefficient (analogous to the tilting module of highest weight $n$). We use the same notation when $\TL(0)$ is defined in characteristic $p$.

\begin{remark} Let us give an independent argument that $s_1$ is a singleton cell. Any closed diagram (e.g. any 1-manifold with empty top and bottom boundary) evaluates to zero in $\TL(0)$.
Consequently, the object $0$ is not a direct summand of the object $n$ for any $n \ge 1$, because the composition of any map $0 \to n$ with any map $n \to 0$ is zero. Similarly, $0$ is not a direct
summand of $T_a \otimes T_b$ for any $a, b \ge 0$ when $a+b \ge 1$. Transferring this observation to the Hecke category via deformation retract, $B_{s_1}$ is not a direct summand of $B_x
\otimes B_y$ for any $x, y$ with length $\ge 1$ where at least one has length $\ge 2$. \end{remark}

\begin{remark} Similarly, let $m=2p^k$ for some $k$, and consider the $\mathbb{Z}^2$-graded Hecke category of the corresponding dihedral group with respect to a $p$-adapted realization (see Corollary \ref{corollary: color graded Hecke category}). We call the basis of the Hecke algebra (with unequal parameters) corresponding to the indecomposable objects the \emph{double-0-$p$ canonical basis}. By the same argument as the previous remark, one can deduce that $\{s_1\}$ and $\{s_2\}$ are singleton cells. Thus the double-0-$p$ canonical basis has (at least) four singleton cells:  $\{\id\}$, $\{s_1\}$, $\{s_2\}$ and $\{w_0\}$, producing all four one-dimensional representations. \end{remark}

\bibliographystyle{plain}
\bibliography{mastercopy}

\end{document}

%% file: preamble.tex
\documentclass[11pt]{amsart}
\usepackage[hmargin=3cm,vmargin=3cm]{geometry}
\usepackage{combelow} 

\usepackage[latin1,utf8]{inputenc} 
\usepackage[english]{babel}
\usepackage{etex}

\usepackage{amssymb,amsfonts,amsthm,amsmath,amscd,latexsym}
\usepackage{palatino, euscript}
\usepackage{graphicx}
\usepackage[all]{xy}
\SelectTips{cm}{}

\usepackage{tikz, tikz-cd}
\usetikzlibrary{matrix, arrows, calc}
\tikzset{frontline/.style={preaction={draw=white,-,line width=6pt}},}  

\usepackage[dvips]{epsfig}
\usepackage{pinlabel}

\newcommand{\ig}[2]{\vcenter{\xy (0,0)*{\includegraphics[scale=#1]{fig/#2}} \endxy}}

\usepackage{enumitem}

\newtheorem{thm}{Theorem}[section]
\newtheorem{lemma}[thm]{Lemma}

\newtheorem{prop}[thm]{Proposition}
\newtheorem{cor}[thm]{Corollary}

\newtheorem{theorem}[thm]{Theorem}

\newtheorem{proposition}[thm]{Proposition}
\newtheorem{corollary}[thm]{Corollary}

\newtheorem*{prop*}{Proposition}

\theoremstyle{definition}
\newtheorem{defn}[thm]{Definition}

\newtheorem{notation}[thm]{Notation}
\newtheorem{ex}[thm]{Example}
\newtheorem{example}[thm]{Example}

\newtheorem{assumption}[thm]{Assumption}
\newtheorem{definition}[thm]{Definition}

\theoremstyle{remark}
\newtheorem{remark}[thm]{Remark}
\newtheorem{rem}[thm]{Remark}
\newtheorem{question}[thm]{Question}

\numberwithin{equation}{section}

\def\CB{{\mathbf C}}  
\def\CC{{\mathcal{C}}}
    
  \def\eb{{\mathbf e}}

    \def\HC{{\mathcal{H}}}

\def\VB{{\mathbf V}}

\let\phi=\varphi

\usepackage{bbm}

\def\Z{{\mathbbm Z}}
\def\Q{{\mathbbm Q}}
\def\one{\mathbbm{1}}

\newcommand{\un}{\underline}

\newcommand{\ot}{\otimes}
\newcommand{\pa}{\partial}
\newcommand{\co}{\colon}

\renewcommand{\to}{\rightarrow}

\newcommand\simto{\xrightarrow{
   \,\smash{\raisebox{-0.65ex}{\ensuremath{\scriptstyle\sim}}}\,}}

\newcommand{\refequal}[1]{\xy {\ar@{=}^{#1}
(-1,0)*{};(1,0)*{}};
\endxy}

\newcommand{\Kar}{\mathbf{Kar}}

\DeclareMathOperator{\ICs}{IC}

\DeclareMathOperator{\Hom}{Hom}

\DeclareMathOperator{\End}{End}

\DeclareMathOperator{\Sym}{Sym}

\DeclareMathOperator{\For}{For}

\DeclareMathOperator{\id}{id}

\DeclareMathOperator{\uni}{uni}
\DeclareMathOperator{\TL}{TL}
\DeclareMathOperator{\JW}{JW}
\DeclareMathOperator{\coeffJW}{coeff}
\DeclareMathOperator{\HeckX}
{\vcenter{\hbox{\begin{tikzpicture}[baseline={(current bounding box.center)}, scale=0.07]
    \draw [dashed] (-6,5) -- (6,5);
    \draw [dashed] (-6,-5) -- (6,-5);
    \draw[red, line width=0.05cm] (-5,-3) -- (-5,-5);
    \draw[red, line width=0.05cm] (-5,3) -- (-5,5);
    \draw[blue, line width=0.05cm] (-3,-3) -- (-3,-5);
    \draw[blue, line width=0.05cm] (-3,3) -- (-3,5);
    \draw[red, line width=0.05cm] (-1,-3) -- (-1,-5);
    \draw[red, line width=0.05cm] (-1,3) -- (-1,5);
    \draw[blue, line width=0.05cm] (5,-3) -- (5,-5);
    \draw[blue, line width=0.05cm] (5,3) -- (5,5);
    \draw (-6,3) rectangle (6,-3);
    \node at (0,0) {X} ;
    \foreach \x in {1,2,3} {
        \node[circle, fill, inner sep=0.01cm] at (\x,-4) {};
        \node[circle, fill, inner sep=0.01cm] at (\x,4) {};
    }
\end{tikzpicture}}}}

\renewcommand{\1}{\mathbbm{1}}

\DeclareMathOperator{\uneq}{uneq}

\DeclareMathOperator{\defect}{defect}
\DeclareMathOperator{\grade}{pos}

\newcommand{\startdotred}{\ig{.5}{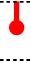}}

\newcommand{\finaldotred}{\ig{.5}{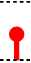}}

\newcommand{\splitred}{\ig{.5}{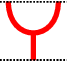}}
\newcommand{\mergered}{\ig{.5}{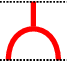}}
\newcommand{\pitchred}{\ig{.5}{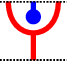}}
\newcommand{\pitchdownred}{\ig{.5}{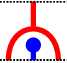}}
\newcommand{\pitchblue}{\ig{.5}{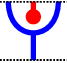}}
\newcommand{\pitchdownblue}{\ig{.5}{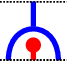}}

\newcommand{\sbotttop}{{
\labellist
\tiny\hair 2pt
\pinlabel {$\dots$} [ ] at 20 27
\pinlabel {$\dots$} [ ] at 20 4
\endlabellist
\centering
\ig{1}{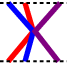}
}}